\documentclass[11pt]{article}
\usepackage[T1]{fontenc}
\usepackage[english]{babel}
\usepackage{microtype}
\usepackage{amsmath,amssymb,amsfonts,amsthm,mathtools}
\usepackage[shortlabels]{enumitem}
\usepackage{geometry}
\usepackage{xcolor}
\usepackage[colorlinks=true,linkcolor=blue,citecolor=blue,urlcolor=blue]{hyperref}
\newtheorem{theorem}{Theorem}[section]
\newtheorem{lemma}[theorem]{Lemma}
\newtheorem{corollary}[theorem]{Corollary}
\newtheorem{proposition}[theorem]{Proposition}

\newtheorem{problem}{Problem}[section]
\newtheorem{claim}{Claim}

\title{Uniform asymptotics for the Kleitman--West problem}
\author{Jiasheng Zeng\thanks{Department of Mathematics, Hong Kong University of Science and Technology, Clear Water Bay, Kowloon, Hong Kong. Email: \href{mailto:jzengbl@connect.ust.hk}{\texttt{jzengbl@connect.ust.hk}}}}
\date{}

\begin{document}
\maketitle

\begin{abstract}
Two distinct $r$-element sets are adjacent if their intersection has size $r-1$. The Kleitman--West problem, equivalently formulated as the edge-isoperimetric problem on the Johnson graph $J(n,r)$, asks for the maximum number $M_r(n,m)$ of unordered adjacent pairs in a family of $m$ $r$-element subsets of an $n$-element set, and for the families attaining this maximum. For every fixed integer $r\ge2$, we determine this maximum asymptotically as $n\to\infty$, uniformly over the entire range of family sizes.

Specifically, for each integer $1\le j\le r$, define $F_j(p)=\max\{p^{1+1/j},\,2p-1+(1-p)^{1+1/j}\}$ for $0\le p\le1$, and $G_j(x)=\lfloor x\rfloor+F_j(x-\lfloor x\rfloor)$ for $x\ge0$. With
\[
P_r(n,m)=\sum_{j=1}^r\frac{n^{j+1}}{2(j-1)!}G_j\!\left(\frac{j!m}{n^j}\right)-\frac{r(r-1)}4mn-\frac m2,
\]
we prove that
\[
\lim_{n\to\infty}\max_{2\le m\le\binom nr}\left|\frac{M_r(n,m)}{P_r(n,m)}-1\right|=0.
\]
The remaining cases satisfy $M_r(n,0)=M_r(n,1)=0$. The proof first establishes the asymptotic formula for dense families using a sharp inequality for bounded random variables, with equality attained by a random variable taking at most two values. To extend the result uniformly to all family sizes, we extract groups of sets sharing a fixed subset, apply the dense estimates after deleting that subset, and control the adjacent pairs not contained in any one group.
\end{abstract}



\section{Introduction}\label{sec:introduction}

Two distinct $r$-element sets are adjacent if their intersection has size $r-1$. For integers $n\ge r\ge1$ and a family $\mathcal F\subseteq\binom{[n]}r$, where $[n]=\{1,\ldots,n\}$, let $A(\mathcal F)$ denote the number of unordered adjacent pairs in $\mathcal F$. For each integer $0\le m\le\binom nr$, define
\begin{equation}\label{eq:extremal}
M_r(n,m)=\max\left\{A(\mathcal F):\mathcal F\subseteq\binom{[n]}r,\ |\mathcal F|=m\right\}.
\end{equation}
Thus $n$ is the size of the ground set, $r$ is the size of each member, and $m$ is the number of members. {The graph case, namely when $r=2$, was solved by Ahlswede and Katona~\cite{AhlswedeKatona1978} in 1978. In the same paper, they explicitly asked for the maximum number of pairs of triples sharing two elements, thereby formulating the case $r=3$. The general problem is discussed by Kleitman~\cite[pp.~60--61]{Kleitman1979} and is stated as an edge-isoperimetric problem by Ahlswede and Cai~\cite{AhlswedeCai1999}. Its exact formulation asks for both the maximum and the families attaining it.}

\begin{problem}[Kleitman--West]\label{prob:KW}
For every pair of integers $n\ge r\ge1$ and every integer $0\le m\le\binom nr$, determine $M_r(n,m)$ and characterize the families $\mathcal F\subseteq\binom{[n]}r$ of size $m$ that attain this maximum.
\end{problem}

The vertices of the Johnson graph $J(n,r)$ are the members of $\binom{[n]}r$, and two vertices are joined when the corresponding sets are adjacent. Therefore $M_r(n,m)$ is the largest number of edges induced by $m$ vertices of $J(n,r)$. Write $\partial_E\mathcal F$ for the set of edges with one endpoint in $\mathcal F$ and the other outside $\mathcal F$. Each vertex has degree $r(n-r)$, so summing degrees over $\mathcal F$ counts each induced edge twice and each boundary edge once. Consequently,
\[
|\partial_E\mathcal F|=r(n-r)|\mathcal F|-2A(\mathcal F).
\]
For fixed $|\mathcal F|=m$, maximizing $A(\mathcal F)$ is thus equivalent to minimizing $|\partial_E\mathcal F|$. We study this edge-isoperimetric problem by estimating $M_r(n,m)$ as $n\to\infty$ for each fixed integer $r\ge2$, allowing $m$ to vary over the full range $0\le m\le\binom nr$.

\subsection{Historical background}\label{subsec:historical-background}

{When $r=2$, the family $\mathcal F$ is the edge set of a simple graph $G$ on $[n]$. For graphs $H$ and $G$, let $N(H,G)$ denote the number of unlabelled, not necessarily induced subgraphs of $G$ isomorphic to $H$, and let $P_k$ denote the path on $k$ vertices. Each pair of adjacent edges forms a copy of $P_3$. Thus, writing $d_G(v)$ for the degree of $v$, we have
\[
A(\mathcal F)=N(P_3,G)=\sum_{v\in[n]}\binom{d_G(v)}2=\frac12\sum_{v\in[n]}d_G(v)^2-m.
\]
Consequently, the graph case is equivalent both to maximizing the number of copies of $P_3$ and to maximizing the sum of squared degrees among graphs with $n$ vertices and $m$ edges.\par}

{Ahlswede and Katona~\cite{AhlswedeKatona1978} proved that a maximizing graph can always be chosen to be quasi-complete or quasi-star. A quasi-complete graph, also known as a quasi-clique, is obtained from a clique by joining at most one additional vertex to some clique vertices and leaving all other vertices isolated. A quasi-star graph is the complement of a quasi-complete graph. These constructions concentrate edges either within a smaller vertex set or around vertices of large degree. Bey~\cite{Bey2006} gave another proof of the graph result and discussed linear-algebraic bounds for the general Johnson-graph problem. Nikiforov~\cite{Nikiforov2007} obtained an asymptotically sharp upper bound for the degree-square extremal function. For each number $n$ of vertices and number $m$ of edges, {\'A}brego, Fern\'andez-Merchant, Neubauer and Watkins~\cite{AbregoFernandezMerchantNeubauerWatkins2009} determined which of the quasi-complete and quasi-star graphs maximizes the number of adjacent edge pairs, when both are optimal, and which other graphs also attain the maximum. Kuzmanovski and Radcliffe~\cite{KR2024} subsequently classified the optimal shifted graphs through their reflect-push method.\par}

{The formulation in terms of $P_3$ places the graph case within the broader problem of maximizing copies of a fixed graph when both the number of vertices and the number of edges are given. At a fixed positive edge density, Nagy~\cite{Nagy2017} determined the asymptotic maximum for $P_5$, with either the quasi-star or the quasi-clique attaining the leading term according to the density. More recently, Yang and Zeng~\cite{YangZeng2026} established the corresponding result for every fixed path $P_{2k+1}$ with $k\ge1$, including an explicit density at which the two constructions exchange optimality. The same two constructions do not suffice for arbitrary fixed graphs. Day and Sarkar~\cite{DaySarkar2021} disproved the conjecture that a quasi-star or a quasi-clique always asymptotically maximizes the subgraph count. A more general reduction was established by Blekherman and Patel~\cite{BlekhermanPatel2024}, who showed that, at a fixed edge density, the asymptotic maximization problem can be restricted to threshold graphs. These are graphs obtained by successively adding a vertex adjacent either to none or to all of the previously added vertices. Thus the results for odd-order paths identify the optimum within a setting in which more complicated constructions can be necessary.\par}

For higher uniformity, the natural counterparts of the two graph constructions are all $r$-sets inside a smaller ground set and all $r$-sets meeting a specified subset. However, Ahlswede and Cai~\cite{AhlswedeCai1999} disproved Kleitman's conjecture that a recursive construction always attains the maximum, even for families of triples. Their counterexample rules out that proposed construction but does not determine the maximum in Problem~\ref{prob:KW}.

Harper~\cite{Harper1991} studied the problem by formulating a related continuous optimization problem. Harper~\cite{Harper2004} later explained an obstruction to the general argument and formulated a stronger conjecture about the possible local optima. The early preprint of Ahlswede and Cai~\cite{AC1993} and their subsequent paper~\cite{AhlswedeCai2006} develop the continuous approach, including the three-dimensional continuous bound and discrete conclusions for particular parameter choices.

{A related development concerns the number of disjoint pairs in a uniform family. For $n\ge2r$, Erd\H{o}s, Ko and Rado~\cite{EKR1961} proved that a pairwise-intersecting family of $r$-subsets of $[n]$ has at most $\binom{n-1}{r-1}$ members, with equality attained by all $r$-sets containing a fixed element. One may then ask how many disjoint pairs are forced when the family is larger. Bollob\'as and Leader~\cite{BollobasLeader2003} studied this minimization problem and conjectured a description of the extremal families. Katona, Katona and Katona~\cite{KatonaKatonaKatona2012} obtained the sharp disjoint-pair bound for a family with one more member than the intersecting-family threshold. Das, Gan and Sudakov~\cite{DasGanSudakov2016} determined the minimum for wider ranges of family sizes. Balogh, Das, Liu, Sharifzadeh and Tran~\cite{BaloghDasLiuSharifzadehTran2019} subsequently extended these results and exhibited counterexamples to the general conjecture of Bollob\'as and Leader.\par}

{The connection with the Kleitman--West problem becomes exact when a larger intersection is required. For $1\le t\le r-1$, a family is $t$-intersecting if every two distinct members have at least $t$ common elements. Das, Gan and Sudakov~\cite{DasGanSudakov2016} studied the associated Erd\H{o}s--Rademacher problem of minimizing the number of pairs whose intersection has size less than $t$. For $t=r-1$, distinct $r$-sets have intersection size at most $r-1$, and hence
\[
A(\mathcal F)=\binom m2-\left|\left\{\{F,F'\}\in\binom{\mathcal F}{2}:|F\cap F'|<r-1\right\}\right|.
\]
Thus their minimization problem is precisely the Kleitman--West maximization problem in this case. For every fixed constant $C>0$, their results determine $M_r(n,m)$ throughout the range $m\le Cn$, provided that $n$ is sufficiently large in terms of $r$ and $C$.\par}

More generally, Bey~\cite{Bey2003} established an upper bound for sums of squared degrees of subsets in uniform hypergraphs, which in particular bounds the number of adjacent pairs. Another source of estimates is the Boolean cube, since two distinct $r$-sets are adjacent exactly when their incidence vectors differ in two coordinates. Kirshner and Samorodnitsky~\cite{KirshnerSamorodnitsky2021} obtained bounds for fixed-distance pair counts in the cube, while Rashtchian and Raynaud~\cite{RashtchianRaynaud2022} studied edge-isoperimetric inequalities for its powers and derived estimates for the Kleitman--West problem. Such estimates and the exact results in specified ranges do not by themselves provide a sharp relative approximation throughout the full range of $m$.

\subsection{Main results}\label{subsec:main-results}

We seek an approximation to $M_r(n,m)$ whose error is negligible compared with $M_r(n,m)$ itself, uniformly over $2\le m\le\binom nr$. When $m$ has order $n^r$, the extremal count has order $n^{r+1}$, so an error bound of order $n^r$ is sufficient for this purpose. For smaller families, however, the same bound may be comparable with or larger than the extremal count, and a more precise estimate is needed. Moreover, constructions at sizes of order $n^j$ can have parts in which all members contain a common $(r-j)$-set. Deleting that common set leaves a $j$-uniform family, so dense problems in lower dimensions enter the answer. We call a range $m\asymp n^j$, with fixed $1\le j\le r$, a critical scale; here $\asymp$ means that $m/n^j$ is bounded above and below by positive constants. A formula valid at every family size must account for these scales and the sizes between them.

To describe the dense contributions, for each integer $j\ge1$ and $0\le p\le1$, we define
\begin{equation}\label{eq:F}
F_j(p)=\max\left\{p^{1+1/j},\ 2p-1+(1-p)^{1+1/j}\right\}.
\end{equation}
The two expressions are the normalized leading counts of, respectively, a complete $j$-uniform family on a smaller ground set and a family of all $j$-sets meeting a specified subset. Their role as extremal leading terms is established below. In particular, $F_1(p)=p^2$. For $x\ge0$, let $\lfloor x\rfloor$ denote the greatest integer not exceeding $x$, and write $\{x\}=x-\lfloor x\rfloor$ for its fractional part. To combine complete parts with a remaining part of intermediate density, we define
\begin{equation}\label{eq:G}
G_j(x)=\lfloor x\rfloor+F_j(\{x\}).
\end{equation}
For $j<r$, each part is obtained by adjoining a fixed $(r-j)$-set to every member of a $j$-uniform family. The term $\lfloor x\rfloor$ counts complete parts, while $\{x\}$ specifies the density of one remaining part. Different parts have different fixed cores but use the same free ground set, which is disjoint from the union of the cores. Thus each part has access to the full free ground set, rather than a separate large block. Since $F_j(0)=0$ and $F_j(1)=1$, the function $G_j$ is continuous, including at the integers.

Using the functions $G_1,\ldots,G_r$, we define
\begin{equation}\label{eq:P}
P_r(n,m)=\sum_{j=1}^r\frac{n^{j+1}}{2(j-1)!}G_j\!\left(\frac{j!m}{n^j}\right)-\frac{r(r-1)}4mn-\frac m2.
\end{equation}
The term proportional to $mn$ cancels linear contributions counted in more than one dimension. The correction $-m/2$ accounts for counting pairs of distinct members and remains necessary when $m$ is bounded. We prove that $P_r(n,m)>0$ whenever $n\ge r\ge2$ and $2\le m\le\binom nr$, so the ratio in the following theorem is well defined. We now state our main result.

\begin{theorem}\label{thm:main}
For every fixed integer $r\ge2$,
\[
\lim_{n\to\infty}\ \max_{2\le m\le\binom nr}\left|\frac{M_r(n,m)}{P_r(n,m)}-1\right|=0.
\]
\end{theorem}

The remaining sizes satisfy $M_r(n,0)=M_r(n,1)=0$. Thus, for fixed uniformity, Theorem~\ref{thm:main} determines the asymptotic value for Problem~\ref{prob:KW}, with relative error tending to zero uniformly over all nontrivial family sizes.

At a fixed positive density, if $m/\binom nr\to p\in(0,1]$, then
\begin{equation}\label{eq:intro-dense}
M_r(n,m)=\left(\frac{F_r(p)}{2(r-1)!}+o(1)\right)n^{r+1}.
\end{equation}
The larger branch of $F_r(p)$ selects either the complete construction on a smaller ground set or the family of all $r$-sets meeting a specified subset. Adding or removing a lower-order number of members makes the size exactly $m$ without changing the leading count.

At a lower critical scale, if $m/n^j\to c>0$ for some $1\le j<r$, then
\begin{equation}\label{eq:intro-critical}
M_r(n,m)=\left(\frac{G_j(j!c)}{2(j-1)!}+o(1)\right)n^{j+1}.
\end{equation}
Writing $j!c=q+p$ with $q=\lfloor j!c\rfloor$ and $0\le p<1$, a construction attaining this leading term has $q$ complete $j$-uniform parts and one part of density $p$, each extended by a different fixed $(r-j)$-set. All parts use a common free ground set of size $n-O_{r,j,c}(1)$. Between consecutive critical scales, if $n^{j-1}=o(m)$ and $m=o(n^j)$ for some $2\le j\le r$, the formula simplifies to $M_r(n,m)=\bigl((j-1)/2+o(1)\bigr)mn$. At the smallest sizes, $M_r(n,m)=\binom m2$ whenever $0\le m\le n-r+1$, since one may take $m$ distinct extensions of the same $(r-1)$-set, making every pair adjacent.

The uniform formula also determines two boundary parameters of the Johnson graph. Put $N=\binom nr$ and let $b_r(n,m)=r(n-r)m-2M_r(n,m)$ be the smallest edge boundary among families of size $m$. The edge-isoperimetric number, without division by the vertex degree, is
\[
h_E(J(n,r))=\min_{1\le m\le\lfloor N/2\rfloor}\frac{b_r(n,m)}m.
\]
For every fixed $r\ge2$, Theorem~\ref{thm:main} implies
\begin{equation}\label{eq:intro-isoperimetric}
h_E(J(n,r))=\left(r\bigl(1-2^{-1/r}\bigr)+o(1)\right)n.
\end{equation}
When $r=2$, this becomes $h_E(J(n,2))=(2-\sqrt2+o(1))n$, recovering the asymptotic result of Fabila-Monroy and Gregorio-Longino~\cite{FG2025}. In this direction, Abiad, van de Berg, Juliano, Reijnders, Simoens, van Veluw and Wittebol~\cite{Abiad2027} study edge-isoperimetric numbers of general graphs and their powers, including results for $J(n,2)$.

To justify that the minimum in~\eqref{eq:intro-isoperimetric} includes sparse sizes, observe that $F_j$ is convex, being the maximum of two convex functions, and has endpoint values $F_j(0)=0$ and $F_j(1)=1$. Thus $F_j(u)\le u$ for $0\le u\le1$ and $G_j(x)\le x$ for $x\ge0$. In~\eqref{eq:P}, bounding the terms with $j<r$ by $jmn/2$ cancels the linear correction, since $\sum_{j=1}^{r-1}j/2=r(r-1)/4$. For $p=r!m/n^r$ and $2\le m\le\lfloor N/2\rfloor$, we therefore have
\[
P_r(n,m)\le\frac{n^{r+1}}{2(r-1)!}F_r(p)-\frac m2
\le\frac{r mn}{2}\,2^{-1/r}.
\]
For the last inequality, $p\le1/2$, and convexity of $F_r$ together with $F_r(0)=0$ implies $F_r(p)\le2pF_r(1/2)=2^{-1/r}p$. The uniform relative estimate in Theorem~\ref{thm:main} now shows that $b_r(n,m)/m\ge\bigl(r(1-2^{-1/r})-o(1)\bigr)n$ uniformly over these sizes. The omitted size $m=1$ has $b_r(n,1)=r(n-r)$. Taking $m=\lfloor N/2\rfloor$ in~\eqref{eq:intro-dense} gives the matching upper bound and proves~\eqref{eq:intro-isoperimetric}.

The bisection width $\operatorname{bw}(J(n,r))$ is the minimum number of edges between two parts whose sizes differ by at most one. Consequently it equals $b_r(n,\lfloor N/2\rfloor)$. Since $F_r(1/2)=2^{-1-1/r}$, the same substitution into~\eqref{eq:intro-dense} yields
\begin{equation}\label{eq:intro-bisection}
\operatorname{bw}(J(n,r))=\left(\frac{1-2^{-1/r}}{2(r-1)!}+o(1)\right)n^{r+1}.
\end{equation}
Unlike the minimum defining $h_E$, this quantity fixes the family size at half of the Johnson graph.

The dense proof begins with the shifting technique surveyed by Frankl~\cite{Frankl1987}. Compression preserves the family size and does not decrease the number of adjacent pairs. For a shifted family, the adjacency count becomes a sum of element labels. Following the connection between compression and continuous optimization developed by Harper~\cite{Harper1991,Harper2004}, we replace this sum by a coordinate integral over a subset of $[0,1]^r$ that is invariant under coordinate permutations and closed under decreasing coordinates. The approximation has an error uniform in the family size. Slicing this set at a uniformly chosen coordinate value produces a random slice volume in $[0,1]$. A sharp one-dimensional probability inequality controls the resulting integrals and makes induction on the dimension possible. Lemma~\ref{lem:general-moment}, proved in Section~\ref{sec:moment}, states the general inequality; its specialization establishes the continuous upper bound used in~\eqref{eq:intro-dense}.

At a lower critical scale $m\asymp n^j$, the sparse proof instead finds common $(r-j)$-sets and groups the members containing them. For fixed approximation parameters, only boundedly many such parts are needed. Removing a common set reduces each part to a $j$-uniform family, to which the dense bound applies. For $j\ge2$, a counting inequality based on fractional edge covers, as in Friedgut and Kahn~\cite{FriedgutKahn1998}, controls the exceptional remaining members containing at least $j$ subsets of size $r-1$ that each occur in many members of the original family. The rest of the remaining family need not be small; its adjacency count is bounded by a linear term, which is combined with the dense contributions by convexity. At the linear scale $j=1$, the necessary estimate follows directly by extracting common $(r-1)$-sets and counting their extensions. A separate upper bound treats sizes between consecutive powers of $n$. For the matching lower bounds, we use one construction in which, for a suitable integer $1\le d<r$, the subfamilies have different fixed $(r-d)$-sets and the same remaining ground set. After removing the fixed set, each subfamily consists of all $d$-sets on that ground set, apart from one subfamily chosen according to the two expressions defining $F_d$. The number of complete subfamilies is fixed at a critical scale and may tend to infinity between consecutive scales. A subsequence argument then combines all scales into the uniform conclusion of Theorem~\ref{thm:main}.

\subsection{Other related problems}\label{subsec:related-problems}

The relation with degree-square problems extends beyond the graph case. For any $S\subseteq[n]$, define $d_{\mathcal F}(S)=|\{F\in\mathcal F:S\subseteq F\}|$. When $|S|=r-1$, this number is the codegree of $S$, and
\[
\sum_{S\in\binom{[n]}{r-1}}d_{\mathcal F}(S)^2=2A(\mathcal F)+rm.
\]
Indeed, each adjacent pair has a unique intersection of size $r-1$, whereas each member contains exactly $r$ subsets of that size. Thus, for fixed $m$, the Kleitman--West problem is exactly the maximization of this codegree squared sum.

The same sum is also studied under conditions that forbid specified subhypergraphs. Balogh, Clemen and Lidick\'y~\cite{BCLSurvey2022} survey these extremal problems. In particular, Balogh, Clemen and Lidick\'y~\cite{BaloghClemenLidicky2022} determined the asymptotic maximum for triple systems containing no complete triple system on four vertices, and likewise on five vertices. Brooks and Linz~\cite{BrooksLinz2026} obtained further exact and asymptotic results under intersection, matching, path and cycle restrictions.

Related results also concern intersections and higher powers of degrees. Wu and Zhang~\cite{WuZhang2026} established a sharp bound for the codegree squared sum of $t$-intersecting families in the range $n\ge(t+1)(r-t+1)$. Cao, Lu and Zhang~\cite{CaoLuZhang2026} extended degree-power estimates to real exponents under intersection and matching restrictions. Zhao and Wang~\cite{ZhaoWang2026} considered sums of powers of degrees of $t$-subsets under intersection, matching and forbidden-path conditions. Chen, I\v{l}kovi\v{c}, Le\'on, Liu and Pikhurko~\cite{CILLP2024} developed a theory of such degree-power sums for nondegenerate Tur\'an problems, in which a fixed subhypergraph is forbidden and the extremal number of edges has positive limiting density. Bao, Chen and Zhang~\cite{BCZ2026} studied analogous squared-codegree questions for intersecting families in finite vector spaces and related finite spaces. These results provide a broader setting for the problem considered here.

Another related problem asks for the maximum sum of pairwise intersection sizes among families of a given size. Huang, Katona and Yue~\cite{HKY2026} studied this quantity and established the optimality of lexicographic constructions in specified parameter ranges. Writing $d_{\mathcal F}(v)=d_{\mathcal F}(\{v\})$, one has
\[
\sum_{\{F,F'\}\in\binom{\mathcal F}{2}}|F\cap F'|
=\sum_{v\in[n]}\binom{d_{\mathcal F}(v)}2.
\]
The right-hand side counts each pair once for each common vertex. Since $\sum_{v\in[n]}d_{\mathcal F}(v)=rm$, maximizing the sum of pairwise intersection sizes is equivalent to maximizing $\sum_{v\in[n]}d_{\mathcal F}(v)^2$ when $m$ is fixed. The sum of pairwise intersection sizes equals $A(\mathcal F)$ when $r=2$, but generally differs from $A(\mathcal F)$ when $r\ge3$.

\vspace{\baselineskip}

The remainder of the paper is organized as follows. In Section~\ref{sec:moment}, we establish a sharp inequality for bounded random variables. In Section~\ref{sec:dense}, we combine this inequality with shifting and a continuous approximation to determine the asymptotic maximum for families of positive limiting density. In Section~\ref{sec:scales}, we obtain bounds for smaller families and construct families attaining the resulting asymptotic estimates. We then combine the estimates for all family sizes to complete the proof of Theorem~\ref{thm:main}.

\section{A sharp inequality for bounded random variables}\label{sec:moment}

For a continuous function $h:[0,1]\to\mathbb R$, we seek the minimum of $\mathbb Eh(X)$ over random variables $X$ taking values in $[0,1]$, when the mean and the expected absolute difference between two independent copies are fixed. For such an $X$, let $X'$ be an independent copy and write $\Gamma(X)=\tfrac12\mathbb E|X-X'|$. If $\mu$ is the probability law of $X$, we also use $\Gamma(\mu)=\tfrac12\iint|x-y|\,d\mu(x)\,d\mu(y)$. The identity $\min(X,X')=(X+X'-|X-X'|)/2$ implies $\mathbb E\min(X,X')=\mathbb EX-\Gamma(X)$. If $p=\mathbb EX$, the inequality $|x-y|\le x+y-2xy$ for $x,y\in[0,1]$, together with independence, implies
\begin{equation}\label{eq:gini-range}
0\le\Gamma(X)\le p(1-p).
\end{equation}
A function $f:[0,1]\to\mathbb R$ is convex if, for every $x,y\in[0,1]$ and every $\lambda\in[0,1]$, we have $f(\lambda x+(1-\lambda)y)\le\lambda f(x)+(1-\lambda)f(y)$. It is concave if the reverse inequality holds for every such choice of $x$, $y$ and $\lambda$. It is strictly concave if, for every pair of distinct points $x,y\in[0,1]$ and every $\lambda\in(0,1)$, we have $f(\lambda x+(1-\lambda)y)>\lambda f(x)+(1-\lambda)f(y)$.

Ba\'illo, C\'arcamo and Mora-Corral~\cite{BailloCarcamoMoraCorral2022} studied related extremal problems for convex functions with fixed endpoint values and a fixed integral, reducing the optimization to finitely many parameters. Here we establish a sharp lower bound for $\mathbb E h(X)$ for a class of symmetric concave functions $h$, with both $\mathbb EX$ and $\Gamma(X)$ fixed. The attaining laws below depend only on these two prescribed quantities and not on the particular function $h$. For $a\in[0,1]$, let $\delta_a$ denote the probability measure on $[0,1]$ satisfying $\delta_a(\{a\})=1$.

\begin{lemma}\label{lem:general-moment}
Let $h:[0,1]\to\mathbb R$ be continuous and concave, with $h(0)=h(1)=0$ and $h(x)=h(1-x)$ for every $x\in[0,1]$. Suppose that $c(x)=h(x)/\sqrt{x(1-x)}$, defined for $0<x<1$, is nondecreasing on $(0,1/2]$. Let $X$ take values in $[0,1]$, and put $p=\mathbb EX$, $q=\max\{p,1-p\}$, $\gamma=\Gamma(X)$ and $K=p(1-p)-\gamma$. Then
\begin{equation}\label{eq:general-moment-sharp}
\mathbb E h(X)\ge\frac{q-\gamma}{q}\,h\!\left(\frac{K}{q-\gamma}\right).
\end{equation}
For every $p\in[0,1]$ and $\gamma\in[0,p(1-p)]$, this bound is attained. More precisely, set $A=K/(q-\gamma)$ and $t_*=(q-\gamma)/q$. An attaining law is
\begin{equation}\label{eq:general-attaining-law}
\mu_*=
\begin{cases}
t_*\delta_A+(\gamma/q)\delta_1,&p\le1/2,\\
t_*\delta_{1-A}+(\gamma/q)\delta_0,&p>1/2.
\end{cases}
\end{equation}
Here $q-\gamma\ge q^2>0$, $0\le A\le1/2$ and $0\le t_*\le1$, so all parameters are well defined. Zero masses and coincident support points are allowed.
\end{lemma}

\begin{proof}[Proof of Lemma~\ref{lem:general-moment}]
Fix an arbitrary random variable $X$ taking values in $[0,1]$, let $\mu_X$ be its law, and let $X'$ be an independent copy. Retain $p=\mathbb EX$, $q=\max\{p,1-p\}$, $\gamma=\tfrac12\mathbb E|X-X'|$ and $K=p(1-p)-\gamma$. We first prove~\eqref{eq:general-moment-sharp} for this $X$. Concavity and the zero endpoint values imply $h(x)\ge(1-x)h(0)+xh(1)=0$ for $0\le x\le1$. Moreover, the pointwise inequality $|x-y|\le x+y-2xy$ and independence imply
\[
0\le\gamma
\le\frac12\mathbb E(X+X'-2XX')
=p-p^2=p(1-p).
\]
Since $p(1-p)=q(1-q)$, it follows that $q-\gamma\ge q^2>0$ and
\[
0\le\frac{K}{q-\gamma}
=\frac{q(1-q)-\gamma}{q-\gamma}
=1-q-\frac{q\gamma}{q-\gamma}
\le1-q\le\frac12.
\]
Thus the right-hand side of~\eqref{eq:general-moment-sharp} is well defined.

If $p=0$ or $p=1$, then $0\le X\le1$ forces $X=p$ almost surely, so $\gamma=K=0$ and both sides of~\eqref{eq:general-moment-sharp} are zero. If $\gamma=0$, independence and $|X-X'|\le1$ imply
\[
0\le2\operatorname{Var}(X)
=\mathbb E(X-X')^2
\le\mathbb E|X-X'|=0.
\]
Hence $X=p$ almost surely, and symmetry of $h$ shows directly that
\[
\frac{q-\gamma}{q}h\!\left(\frac{K}{q-\gamma}\right)
=h\!\left(\frac{p(1-p)}q\right)
=h(\min\{p,1-p\})=h(p)=\mathbb Eh(X).
\]
If $\gamma=p(1-p)$, then $K=0$, so the right-hand side of~\eqref{eq:general-moment-sharp} is zero and the inequality follows from $h\ge0$. It remains to consider $0<p<1$ and $0<\gamma<p(1-p)$. First suppose that $h$ is strictly concave, so $h(x)>0$ for $0<x<1$.

To obtain a lower bound for $\mathbb Eh(X)$, let $\mathcal P_{p,\gamma}$ be the set of probability measures $\nu$ on $[0,1]$ satisfying $\int x\,d\nu(x)=p$ and $\Gamma(\nu)=\tfrac12\iint_{[0,1]^2}|x-y|\,d\nu(x)\,d\nu(y)\le\gamma$, and put $L=\inf_{\nu\in\mathcal P_{p,\gamma}}\int h\,d\nu$. The law $\mu_X$ satisfies both constraints, and therefore
\[
\mathbb Eh(X)=\int h\,d\mu_X\ge L\ge0.
\]
All distributional reductions below concern auxiliary minimizers for $L$, rather than the original law $\mu_X$.

We first prove that the infimum is attained.

\begin{claim}\label{clm:moment-minimizer}
There exists $\nu\in\mathcal P_{p,\gamma}$ such that $\int h\,d\nu=L$.
\end{claim}

\begin{proof}[Proof of Claim~\ref{clm:moment-minimizer}]
Choose $\nu_n\in\mathcal P_{p,\gamma}$ with $\int h\,d\nu_n\to L$, and define the nondecreasing functions
\[
g_n(u)=\inf\{x\in[0,1]:\nu_n([0,x])\ge u\},
\qquad 0<u<1.
\]
If $U$ is uniform on $(0,1)$, right-continuity of the distribution function implies $\mathbb P(g_n(U)\le x)=\mathbb P(U\le\nu_n([0,x]))=\nu_n([0,x])$ for $0\le x\le1$. Thus $g_n(U)$ has law $\nu_n$.

Enumerate the rational numbers in $(0,1)$ as $u_1,u_2,\ldots$. Since $0\le g_n\le1$, we can successively choose nested subsequences so that the $j$th subsequence converges at $u_1,\ldots,u_j$. Taking the $k$th term of the $k$th subsequence produces a common subsequence converging at every $u_j$, because for each fixed $j$, all terms with $k\ge j$ belong to the $j$th subsequence. Retain the notation $(\nu_n,g_n)$ for the selected subsequence, and denote the limiting value at $u_j$ by $b_j=\lim_{n\to\infty}g_n(u_j)$.

Define $g(u)=\sup\{b_j:u_j<u\}$ for $0<u<1$. Then $0\le g\le1$, and $g$ is nondecreasing. For $u_i<u$, monotonicity implies $b_i\le\liminf_n g_n(u)$. For any $u<v<1$, choose $u_j\in(u,v)$; then $\limsup_n g_n(u)\le b_j\le g(v)$. Taking the supremum over $u_i<u$ and letting $v\downarrow u$, we obtain
\[
g(u)=\sup_{u_i<u}b_i\le\liminf_{n\to\infty}g_n(u)\le\limsup_{n\to\infty}g_n(u)\le g(u+),
\]
where $g(u+)$ is the right-hand limit of $g$ at $u$. Thus $g_n(u)\to g(u)$ at every continuity point of $g$. Since a bounded nondecreasing function has at most countably many discontinuities, $g_n(u)\to g(u)$ for almost every $u\in(0,1)$.

Let $\nu$ be the law of $g(U)$. Since $g_n\to g$ almost everywhere and $0\le g_n\le1$, the dominated convergence theorem allows us to pass to the limit in the mean constraint, obtaining
\[
\int x\,d\nu(x)=\int_0^1g(u)\,du=\lim_{n\to\infty}\int_0^1g_n(u)\,du=p.
\]
For the second constraint, let $E\subseteq(0,1)$ be a set of Lebesgue measure zero such that $g_n(u)\to g(u)$ for every $u\notin E$. Then $|g_n(u)-g_n(v)|\to|g(u)-g(v)|$ whenever $u,v\notin E$. The exceptional pairs therefore lie in $(E\times(0,1))\cup((0,1)\times E)$, which has two-dimensional Lebesgue measure zero. Since $|g_n(u)-g_n(v)|\le1$, the dominated convergence theorem implies
\[
\Gamma(\nu)=\frac12\int_0^1\int_0^1|g(u)-g(v)|\,du\,dv=\lim_{n\to\infty}\Gamma(\nu_n)\le\gamma.
\]
Thus $\nu\in\mathcal P_{p,\gamma}$. Finally, continuity and boundedness of $h$ allow dominated convergence for the objective, so
\[
\int h\,d\nu=\int_0^1h(g(u))\,du=\lim_{n\to\infty}\int_0^1h(g_n(u))\,du=L.\qedhere
\]
\end{proof}

We next show that every minimizer satisfies $\Gamma(\nu)=\gamma$.

\begin{claim}\label{clm:moment-constraint-equality}
If $\nu\in\mathcal P_{p,\gamma}$ satisfies $\int h\,d\nu=L$, then $\Gamma(\nu)=\gamma$.
\end{claim}

\begin{proof}[Proof of Claim~\ref{clm:moment-constraint-equality}]
Suppose, for a contradiction, that some $\nu\in\mathcal P_{p,\gamma}$ satisfies $\int h\,d\nu=L$ and $\Gamma(\nu)<\gamma$. To construct another probability law with mean $p$, consider $\beta_p=(1-p)\delta_0+p\delta_1$, which has mean $p$ and $\Gamma(\beta_p)=p(1-p)$. Since $\beta_p$ assigns probabilities $1-p$ and $p$ to $0$ and $1$, respectively, and $\nu$ has mean $p$, integration first with respect to $\beta_p$ shows that
\[
\begin{aligned}
\iint|x-y|\,d\nu(x)\,d\beta_p(y)
&=\int\bigl((1-p)x+p(1-x)\bigr)\,d\nu(x)\\
&=(1-p)p+p(1-p)=2p(1-p).
\end{aligned}
\]
For $0\le s\le1$, the measure $\nu_s=(1-s)\nu+s\beta_p$ is a probability law with mean $\int x\,d\nu_s(x)=(1-s)p+sp=p$. Expanding the definition of $\Gamma$ and using the preceding integral, we obtain
\begin{equation}\label{eq:moment-mixture}
\begin{aligned}
\Gamma(\nu_s)
&=(1-s)^2\Gamma(\nu)+2s(1-s)p(1-p)+s^2p(1-p)\\
&=p(1-p)-(1-s)^2\bigl(p(1-p)-\Gamma(\nu)\bigr).
\end{aligned}
\end{equation}
Since $\Gamma(\nu)<\gamma<p(1-p)$, we can choose
\[
s=1-\sqrt{\frac{p(1-p)-\gamma}{p(1-p)-\Gamma(\nu)}}\in(0,1).
\]
Substituting this value into~\eqref{eq:moment-mixture} gives $\Gamma(\nu_s)=\gamma$. Together with the mean calculation, this proves $\nu_s\in\mathcal P_{p,\gamma}$.

To compare $\int h\,d\nu_s$ with $L$, first note that $\Gamma(\nu)<p(1-p)=\Gamma(\beta_p)$ implies $\nu\ne\beta_p$. The only probability law concentrated on $\{0,1\}$ with mean $p$ is $\beta_p$, so $\nu((0,1))>0$. Since $h(x)>0$ for $0<x<1$, it follows that $L=\int h\,d\nu>0$. The membership $\nu_s\in\mathcal P_{p,\gamma}$, the zero endpoint values of $h$ and $0<s<1$ therefore imply
\[
L\le\int h\,d\nu_s=(1-s)\int h\,d\nu+s\bigl((1-p)h(0)+ph(1)\bigr)=(1-s)L<L.
\]
This contradicts the definition of $L$ and proves that every minimizer satisfies $\Gamma(\nu)=\gamma$.
\end{proof}

For a probability law $\nu$ on $[0,1]$, a point $x\in[0,1]$ is called a support point if $\nu((x-\delta,x+\delta)\cap[0,1])>0$ for every $\delta>0$. The set of all support points is called the support of $\nu$. We next use Claim~\ref{clm:moment-constraint-equality} to show that every minimizer has at most three support points.

\begin{claim}\label{clm:moment-three-points}
Every $\nu\in\mathcal P_{p,\gamma}$ satisfying $\int h\,d\nu=L$ has at most three support points.
\end{claim}

\begin{proof}[Proof of Claim~\ref{clm:moment-three-points}]
Suppose, for a contradiction, that some $\nu\in\mathcal P_{p,\gamma}$ satisfies $\int h\,d\nu=L$ and has at least four support points. By Claim~\ref{clm:moment-constraint-equality}, we have $\Gamma(\nu)=\gamma$. Choose four distinct support points and four pairwise disjoint Borel neighborhoods $B_1,\ldots,B_4$, one around each point. The definition of a support point ensures that $\nu(B_i)>0$ for every $i$. Put $V_\nu(x)=\int|x-y|\,d\nu(y)$. Consider the following homogeneous linear equations in the four unknowns $a_1,\ldots,a_4$,
\[
\begin{cases}
\displaystyle\sum_{i=1}^4 a_i\nu(B_i)=0,\\[4pt]
\displaystyle\sum_{i=1}^4 a_i\int_{B_i}x\,d\nu(x)=0,\\[4pt]
\displaystyle\sum_{i=1}^4 a_i\int_{B_i}V_\nu(x)\,d\nu(x)=0.
\end{cases}
\]
Since this system has three equations and four unknowns, it has a solution with $a_1,\ldots,a_4$ not all zero. Fix such a solution and define the signed measure $\eta$ by $d\eta=\sum_{i=1}^4a_i\mathbf1_{B_i}\,d\nu$. Equivalently, for every Borel set $A\subseteq[0,1]$,
\[
\eta(A)=\int_A\sum_{i=1}^4a_i\mathbf1_{B_i}(x)\,d\nu(x)=\sum_{i=1}^4a_i\,\nu(A\cap B_i).
\]
By the definition of $\eta$, the three equations can be written as
\begin{equation}\label{eq:moment-perturbation-constraints}
\int d\eta=0,\qquad
\int x\,d\eta(x)=0,\qquad
\int V_\nu\,d\eta=0.
\end{equation}
Since the coefficients are not all zero, choose $j$ with $a_j\ne0$. The sets $B_1,\ldots,B_4$ are pairwise disjoint, so $\eta(B_j)=a_j\nu(B_j)\ne0$. Hence $\eta$ is nonzero. Rescale the coefficients so that $\max_i|a_i|\le1$. Since the sets $B_1,\ldots,B_4$ are pairwise disjoint, the sum $\sum_{i=1}^4a_i\mathbf1_{B_i}(x)$ equals either $0$ or one of the coefficients $a_i$, and hence has absolute value at most $1$. For $0<\varepsilon<1$, we therefore have $1\pm\varepsilon\sum_{i=1}^4a_i\mathbf1_{B_i}(x)\ge1-\varepsilon>0$ for every $x\in[0,1]$. Consequently, for every Borel set $A\subseteq[0,1]$,
\[
(\nu\pm\varepsilon\eta)(A)=\int_A\left(1\pm\varepsilon\sum_{i=1}^4a_i\mathbf1_{B_i}(x)\right)d\nu(x)\ge(1-\varepsilon)\nu(A)\ge0.
\]
Thus $\nu+\varepsilon\eta$ and $\nu-\varepsilon\eta$ are nonnegative measures. The first two identities in~\eqref{eq:moment-perturbation-constraints} imply $(\nu\pm\varepsilon\eta)([0,1])=1\pm\varepsilon\eta([0,1])=1$ and $\int x\,d(\nu\pm\varepsilon\eta)(x)=p\pm\varepsilon\int x\,d\eta(x)=p$. Hence both are probability laws with mean $p$. To compute $\Gamma(\nu\pm\varepsilon\eta)$, we use the identity
\[
|x-y|=\int_0^1
\bigl(\mathbf1_{\{x>s\}}-\mathbf1_{\{y>s\}}\bigr)^2\,ds.
\]
Expanding the square and using $\eta([0,1])=0$ removes the two terms that depend on only one variable, leaving
\begin{equation}\label{eq:negative-kernel}
\iint|x-y|\,d\eta(x)\,d\eta(y) =-2\int_0^1\eta((s,1])^2\,ds =-2\int_0^1\eta([0,s])^2\,ds.
\end{equation}
Write $C=\int_0^1\eta([0,s])^2\,ds$. This number is positive. Otherwise, $\eta([0,s])=0$ for almost every $s\in[0,1]$. Since the function $s\mapsto\eta([0,s])$ is right-continuous, it follows that $\eta([0,s])=0$ for every $s\in[0,1)$. Together with $\eta([0,1])=0$, this would force $\eta=0$, a contradiction. Expanding $\Gamma$ and applying~\eqref{eq:moment-perturbation-constraints} and~\eqref{eq:negative-kernel}, we obtain
\begin{equation}\label{eq:gini-perturbation}
\begin{aligned}
\Gamma(\nu\pm\varepsilon\eta)
&=\Gamma(\nu)\pm\varepsilon\int V_\nu\,d\eta
+\frac{\varepsilon^2}{2}\iint|x-y|\,d\eta(x)\,d\eta(y)\\
&=\gamma-\varepsilon^2C<\gamma.
\end{aligned}
\end{equation}
Thus both perturbed laws belong to $\mathcal P_{p,\gamma}$. Since $\int h\,d\nu=L$, linearity of integration implies
\[
\int h\,d(\nu\pm\varepsilon\eta)=\int h\,d\nu\pm\varepsilon\int h\,d\eta=L\pm\varepsilon\int h\,d\eta.
\]
Hence $\int h\,d(\nu\pm\varepsilon\eta)\le L$ for at least one choice of sign. For this choice of sign, the definition of $L$ implies $L\le\int h\,d(\nu\pm\varepsilon\eta)\le L$. Thus $\nu\pm\varepsilon\eta$ is another minimizer, but $\Gamma(\nu\pm\varepsilon\eta)<\gamma$, contradicting Claim~\ref{clm:moment-constraint-equality}. Hence every minimizer has at most three support points.
\end{proof}

To express a minimizer as a finite sum of the measures $\delta_x$, fix a minimizer $\nu\in\mathcal P_{p,\gamma}$ and let $S$ be its support. By Claim~\ref{clm:moment-three-points}, the set $S$ has at most three elements. We first show that $\nu([0,1]\setminus S)=0$. For each $x\in[0,1]\setminus S$, the definition of support provides $\delta>0$ such that $\nu((x-\delta,x+\delta)\cap[0,1])=0$. Choose rational numbers $a,b$ with $x-\delta<a<x<b<x+\delta$. Then $x\in(a,b)\cap[0,1]\subseteq(x-\delta,x+\delta)\cap[0,1]$, so $\nu((a,b)\cap[0,1])=0$. The sets $(a,b)\cap[0,1]$ with rational $a<b$ and $\nu((a,b)\cap[0,1])=0$ form a countable collection, which we denote by $I_1,I_2,\ldots$. These sets cover $[0,1]\setminus S$, so countable subadditivity implies
\[
0\le\nu([0,1]\setminus S)\le\nu\left(\bigcup_{j=1}^{\infty}I_j\right)\le\sum_{j=1}^{\infty}\nu(I_j)=0.
\]
Consequently $\nu(S)=1$, and in particular $S$ is nonempty.

Write $S=\{x_1,\ldots,x_k\}$, where $0\le x_1<\cdots<x_k\le1$ and $1\le k\le3$, and put $w_i=\nu(\{x_i\})$. Since the points are distinct and $S$ is finite, for each $i$ there is $\delta_i>0$ such that $U_i=(x_i-\delta_i,x_i+\delta_i)\cap[0,1]$ satisfies $U_i\cap S=\{x_i\}$. The definition of a support point ensures that $\nu(U_i)>0$, whereas $U_i\setminus\{x_i\}\subseteq[0,1]\setminus S$ has measure zero. Thus $0<\nu(U_i)=\nu(\{x_i\})=w_i$, and finite additivity implies $\sum_{i=1}^kw_i=\nu(S)=1$. For every Borel set $A\subseteq[0,1]$, the set $A\setminus S$ also has measure zero, and hence
\[
\nu(A)=\nu(A\cap S)+\nu(A\setminus S)=\sum_{\substack{1\le i\le k\\x_i\in A}}w_i=\sum_{i=1}^kw_i\delta_{x_i}(A).
\]
Since this identity holds for every Borel set $A$, it proves that $\nu=\sum_{i=1}^kw_i\delta_{x_i}$ with $w_i>0$ and $\sum_{i=1}^kw_i=1$.

The representation above still allows more than one support point in $(0,1)$. The next claim shows that a minimizer can be chosen with at most one support point in $(0,1)$, so that its support is contained in $\{0,a,1\}$ for some $a\in(0,1)$.

\begin{claim}\label{clm:moment-endpoint-support}
There exist $a\in(0,1)$ and $\widehat\nu\in\mathcal P_{p,\gamma}$ such that $\int h\,d\widehat\nu=L$ and $\widehat\nu(\{0,a,1\})=1$.
\end{claim}

\begin{proof}[Proof of Claim~\ref{clm:moment-endpoint-support}]
We first show that a three-point minimizer cannot have all its support points in $(0,1)$. Suppose, for a contradiction, that a minimizer has the form $\nu=\sum_{i=1}^3w_i\delta_{x_i}$, where $0<x_1<x_2<x_3<1$, each $w_i>0$ and $\sum_iw_i=1$. Define $b_i=\sum_{k<i}w_k-\sum_{k>i}w_k$. As long as the points remain ordered, the constraints have the form
\begin{equation}\label{eq:position-constraints}
\int x\,d\nu(x)=\sum_{i=1}^3w_ix_i=p,
\qquad
\Gamma(\nu)=\sum_{i<j}w_iw_j(x_j-x_i)
=\sum_{i=1}^3w_ib_ix_i=\gamma.
\end{equation}
Choose a nonzero solution $(v_1,v_2,v_3)$ to the two homogeneous linear equations $\sum_iw_iv_i=0$ and $\sum_iw_ib_iv_i=0$. For sufficiently small $\varepsilon>0$, both sets of points $x_i+\varepsilon v_i$ and $x_i-\varepsilon v_i$ remain interior and ordered. The laws $\nu_\pm=\sum_{i=1}^3w_i\delta_{x_i\pm\varepsilon v_i}$ consequently retain mean $p$ and $\Gamma=\gamma$. Since at least one $v_i$ is nonzero and all weights are positive, strict concavity implies
\[
L=\sum_{i=1}^3w_ih(x_i) >\frac12\sum_{i=1}^3w_i \bigl(h(x_i+\varepsilon v_i)+h(x_i-\varepsilon v_i)\bigr) =\frac12\left(\int h\,d\nu_++\int h\,d\nu_-\right)\ge L.
\]
This contradiction shows that every three-point minimizer contains an endpoint.

To reduce the remaining possibilities to support contained in $\{0,a,1\}$, a law already containing both endpoints, or a two-point law containing an endpoint, requires no change. A one-point law has $\Gamma=0$ and is excluded by the present assumptions. Thus the remaining possibilities have support $\{0,a,b\}$, $\{a,b\}$ or $\{a,b,1\}$, where $0<a<b<1$. We first consider the first two cases. A minimizer in either case can be written as
\[
\nu=(1-t)\delta_0+(t-w)\delta_a+w\delta_b,
\qquad 0<w<t\le1.
\]
Here $t$ is the total probability assigned to $a$ and $b$, and $w$ is the probability assigned to $b$. When $t<1$, all three coefficients are positive, so the support is $\{0,a,b\}$. When $t=1$, the coefficient of $\delta_0$ vanishes, so the support is $\{a,b\}$. The remaining case with support $\{a,b,1\}$ will be treated by reflection after these two cases. Since $\nu$ is a minimizer, Claim~\ref{clm:moment-constraint-equality} implies $\Gamma(\nu)=\gamma$. Put $\Delta=b-a$ and $Q=p-\gamma$. Since $\int x\,d\nu(x)=p$, the identity $\Gamma(\nu)=p-\iint\min(x,y)\,d\nu(x)\,d\nu(y)$ shows that $Q=p-\Gamma(\nu)=\iint\min(x,y)\,d\nu(x)\,d\nu(y)$. In this double integral, every pair containing $0$ contributes zero, while the pairs $(a,a),(a,b),(b,a),(b,b)$ have minimum $a,a,a,b$, respectively. Thus the mean and $Q$ can be written as
\begin{equation}\label{eq:endpoint-constraints}
\begin{aligned}
p&=(t-w)a+wb=ta+w\Delta,\\
Q&=a(t-w)^2+2a(t-w)w+bw^2=t^2a+w^2\Delta.
\end{aligned}
\end{equation}

To vary $w,a,b$ while keeping $t$ fixed, put $u=1/t$, $z_0=1/w-u$ and $D=p-uQ$. Then $u\ge1$, $z_0>0$, and
\[
D=ta+w\Delta-\frac1t(t^2a+w^2\Delta)
=w\left(1-\frac wt\right)\Delta>0.
\]
Also $Q=p-\gamma\ge p^2>0$. We shall use these parameters to vary the weight $w$ and adjust $a,b$, while keeping the support in $[0,1]$, the mean equal to $p$ and $\Gamma$ equal to $\gamma$. Our aim is to move $a$ to $0$ or $b$ to $1$ without increasing the integral of $h$. A law satisfying these conditions still belongs to $\mathcal P_{p,\gamma}$ and has integral of $h$ at most $L$. By the definition of $L$, this integral is also at least $L$, so the law is another minimizer and has at most one support point in $(0,1)$.

Keep $u,p,\gamma$, and hence $D,Q$, fixed. Replace the original weight $w$ by $w(z)=1/(u+z)$, where $z>0$. The original weight is recovered at $z=z_0$, since $w(z_0)=1/(u+z_0)=w$, and $0<w(z)<1/u=t$ for every $z>0$. We next determine the new positions $a(z),b(z)$ by requiring the two equalities in~\eqref{eq:endpoint-constraints} to remain valid. Writing $\Delta(z)=b(z)-a(z)$, these equalities become $p=a(z)/u+w(z)\Delta(z)$ and $Q=a(z)/u^2+w(z)^2\Delta(z)$. Subtracting $u$ times the second equality from the first leads to $D=w(z)(1-uw(z))\Delta(z)$. Since $w(z)(1-uw(z))=z/(u+z)^2>0$, we obtain $\Delta(z)=D(u+z)^2/z$. Substituting this expression in the equality for $Q$ determines $a(z)$, and then $b(z)=a(z)+\Delta(z)$. Using $D=p-uQ$ to simplify the latter expression, we find
\begin{equation}\label{eq:endpoint-parametrization}
a(z)=u^2Q-\frac{u^2D}{z},
\qquad
b(z)=a(z)+\Delta(z)=u(2p-uQ)+Dz.
\end{equation}
Since $\Delta(z)>0$, the points belong to $[0,1]$ exactly when $a(z)\ge0$ and $b(z)\le1$. These inequalities are equivalent to
\begin{equation}\label{eq:endpoint-feasible}
z_-:=\frac DQ\le z\le
z_+:=\frac{1-u(2p-uQ)}D.
\end{equation}
At the original value $z_0$, the support points are $a(z_0)=a>0$ and $b(z_0)=b<1$, so $z_-<z_0<z_+$. In particular, the interval is nonempty. For every $z$ in this interval, define
\[
\nu_z=
\left(1-\frac1u\right)\delta_0
+\frac{z}{u(u+z)}\delta_{a(z)}
+\frac1{u+z}\delta_{b(z)}.
\]
Its weights are nonnegative and sum to $1-1/u+(z+u)/(u(u+z))=1$. Its mean is
\[
\int x\,d\nu_z(x) =\frac{a(z)}u+\frac{\Delta(z)}{u+z} =uQ-\frac{uD}{z}+\frac{D(u+z)}z =uQ+D=p.
\]
The same calculation of the minimum as in~\eqref{eq:endpoint-constraints} shows that
\[
\iint\min(x,y)\,d\nu_z(x)\,d\nu_z(y) =\frac{a(z)}{u^2}+\frac{\Delta(z)}{(u+z)^2} =Q-\frac Dz+\frac Dz=Q.
\]
Consequently $\Gamma(\nu_z)=p-Q=\gamma$, so every $\nu_z$ is feasible and $\nu_{z_0}=\nu$.

The objective of $\nu_z$ can be written as a quotient with a positive affine denominator. Since $h(0)=0$, set
\begin{equation}\label{eq:perspective-objective}
J(z)=\int h\,d\nu_z=\frac{N(z)}{u+z},
\qquad
N(z)=\frac zu\,h(a(z))+h(b(z)).
\end{equation}
The numerator $N$ is concave. Indeed, for $z=\lambda z_1+(1-\lambda)z_2$, where $z_1,z_2\in[z_-,z_+]$ and $0\le\lambda\le1$, the expression for $a(z)$ implies
\[
a(z)=\frac{\lambda z_1}{z}a(z_1)
+\frac{(1-\lambda)z_2}{z}a(z_2).
\]
The two coefficients are nonnegative and sum to one. Concavity of $h$ therefore implies
\[
\frac zu\,h(a(z))
\ge\frac{\lambda z_1}{u}h(a(z_1))
+\frac{(1-\lambda)z_2}{u}h(a(z_2)).
\]
Since $b$ is affine, we also have $h(b(z))\ge\lambda h(b(z_1))+(1-\lambda)h(b(z_2))$. Adding these inequalities proves $N(z)\ge\lambda N(z_1)+(1-\lambda)N(z_2)$. Dividing by $u+z=\lambda(u+z_1)+(1-\lambda)(u+z_2)>0$, we conclude that
\[
J(z)\ge\frac{\lambda(u+z_1)J(z_1)+(1-\lambda)(u+z_2)J(z_2)}{\lambda(u+z_1)+(1-\lambda)(u+z_2)}\ge\min\{J(z_1),J(z_2)\}.
\]
Apply this inequality with $z_1=z_-$, $z_2=z_+$ and $z=z_0$. Since the original law minimizes the objective and both endpoint laws are feasible, $L\le\min\{J(z_-),J(z_+)\}\le J(z_0)=L$. Thus at least one endpoint law is also a minimizer. At $z_-$ the lower support point is $a(z_-)=0$, whereas at $z_+$ the upper support point is $b(z_+)=1$. Either choice is supported on $\{0,a,1\}$ for some $a\in(0,1)$, with a zero endpoint weight allowed.

If the original three-point minimizer contains $1$ but not $0$, replace each support point $x$ by $1-x$, retaining its weight. The reflected law has mean $1-p$, the same value of $\Gamma$, and the same objective by symmetry of $h$. It minimizes the objective for its mean, since reflecting any feasible law with a smaller objective back would contradict the minimality of the original law. Apply the preceding reduction to this law and then reflect back. The resulting law has mean $p$ and $\Gamma=\gamma$, is supported on $\{0,a,1\}$, and has objective at most $L$, hence equal to $L$. This covers all remaining support configurations. Therefore there is an auxiliary minimizer $\widehat\nu$ supported on $\{0,a,1\}$.
\end{proof}

To bound $L$, choose $a$ and $\widehat\nu$ as in Claim~\ref{clm:moment-endpoint-support}. Since $\widehat\nu$ is a minimizer, Claim~\ref{clm:moment-constraint-equality} ensures that $\Gamma(\widehat\nu)=\gamma$. First suppose that $p\le1/2$, so $q=1-p$. If $t$ is the probability assigned to $a\in(0,1)$, the mean and total probability determine the endpoint probabilities, and thus
\[
\widehat\nu=\bigl(q-t(1-a)\bigr)\delta_0
+t\delta_a+(p-ta)\delta_1.
\]
Here $t>0$, since a law supported only on the endpoints and having mean $p$ would have $\Gamma=p(1-p)>\gamma$. In the integral of the minimum, a pair containing zero contributes zero, so
\[
\iint\min(x,y)\,d\widehat\nu(x)\,d\widehat\nu(y) =a\bigl(t^2+2t(p-ta)\bigr)+(p-ta)^2 =p^2+t^2a(1-a).
\]
Hence $\gamma=p-p^2-t^2a(1-a)=pq-t^2a(1-a)$, and therefore $K=t^2a(1-a)$. The endpoint values of $h$ now imply
\begin{equation}\label{eq:three-point-objective}
\int h\,d\widehat\nu
=th(a)
=\sqrt{\frac{K}{a(1-a)}}\,h(a)
=\sqrt K\,c(a).
\end{equation}
Nonnegativity of the endpoint probabilities means $ta\le p$ and $t(1-a)\le q$. Squaring and substituting $t^2=K/(a(1-a))$ produces $Ka/(1-a)\le p^2$ and $K(1-a)/a\le q^2$. Equivalently,
\begin{equation}\label{eq:three-point-feasible}
A=\frac K{K+q^2}\le a\le
B:=\frac{p^2}{p^2+K}.
\end{equation}
Here $K+q^2=pq-\gamma+q^2=q-\gamma$, so $A$ is exactly $K/(q-\gamma)$ as in the statement. Since $p\le q$, the upper bound on $a$ also implies
\[
1-a\ge1-B=\frac K{p^2+K}
\ge\frac K{q^2+K}=A.
\]
Thus $0<A\le\min\{a,1-a\}\le1/2$. Symmetry of $h$ implies $c(a)=c(1-a)$, and the assumed monotonicity of $c$ gives
\[
\int h\,d\widehat\nu
=\sqrt K\,c(a)
=\sqrt K\,c(\min\{a,1-a\})
\ge\sqrt K\,c(A).
\]
To evaluate this lower bound, the identity $K+q^2=q-\gamma$ implies
\begin{equation}\label{eq:moment-extremizer}
A=\frac K{q-\gamma},
\qquad
1-A=\frac{q^2}{q-\gamma},
\qquad
A(1-A)=\frac{Kq^2}{(q-\gamma)^2}.
\end{equation}
Since the original law $\mu_X$ is feasible and $\widehat\nu$ minimizes the objective, we obtain the comparison
\[
\begin{aligned}
\mathbb Eh(X)
&=\int h\,d\mu_X
\ge L=\int h\,d\widehat\nu\\
&\ge\sqrt K\,c(A)
=\frac{\sqrt K}{\sqrt{A(1-A)}}\,h(A)
=\frac{q-\gamma}{q}\,
h\!\left(\frac K{q-\gamma}\right).
\end{aligned}
\]
This proves~\eqref{eq:general-moment-sharp} for $p\le1/2$. If $p>1/2$, apply this conclusion to $1-X$. Its mean is $1-p<1/2$, its values of $q,\gamma,K$ are unchanged, and $\mathbb Eh(1-X)=\mathbb Eh(X)$ by symmetry. The same bound therefore holds for the original $X$.

It remains to remove the strict concavity assumption from the inequality. For $\varepsilon>0$, put $h_\varepsilon(x)=h(x)+\varepsilon x(1-x)$. This function is continuous, symmetric and zero at both endpoints. For distinct $x,y\in[0,1]$ and $0<\lambda<1$, concavity of $h$ implies
\[
h_\varepsilon(\lambda x+(1-\lambda)y) -\lambda h_\varepsilon(x)-(1-\lambda)h_\varepsilon(y) \ge\varepsilon\lambda(1-\lambda)(x-y)^2>0.
\]
Thus $h_\varepsilon$ is strictly concave. Moreover, $h_\varepsilon(x)/\sqrt{x(1-x)}=c(x)+\varepsilon\sqrt{x(1-x)}$ for $0<x<1$, and both summands are nondecreasing on $(0,1/2]$. Applying the strictly concave case to $h_\varepsilon$ and the original $X$, with $A=K/(q-\gamma)$, we obtain
\[
\mathbb Eh(X)+\varepsilon\mathbb E[X(1-X)]
\ge\frac{q-\gamma}{q}
\bigl(h(A)+\varepsilon A(1-A)\bigr).
\]
Since $0\le X(1-X)\le1/4$ and $0\le A(1-A)\le1/4$, letting $\varepsilon\downarrow0$ proves~\eqref{eq:general-moment-sharp} for every $h$ in the statement.

Finally, we verify attainment independently of the lower-bound argument. Fix any $p\in[0,1]$ and $\gamma\in[0,p(1-p)]$, retain $q=\max\{p,1-p\}$ and $K=p(1-p)-\gamma$, and set $A=K/(q-\gamma)$ and $t_*=(q-\gamma)/q$. The parameter bounds above imply $0\le A\le1/2$, $0\le t_*\le1$ and $t_*+\gamma/q=1$, so the expressions in~\eqref{eq:general-attaining-law} define probability laws. If $p\le1/2$, let $Y$ be a new random variable with law $\mu_*=t_*\delta_A+(\gamma/q)\delta_1$. Since $q=1-p$, its mean is $\mathbb EY=t_*A+\gamma/q=(K+\gamma)/q=p$. Also $1-A=q^2/(q-\gamma)$, so $\Gamma(Y)=t_*(\gamma/q)(1-A)=\gamma$. Finally, $h(1)=0$ implies $\mathbb Eh(Y)=t_*h(A)=((q-\gamma)/q)h(K/(q-\gamma))$, as required.

If $p>1/2$, instead take $Y$ with law $\mu_*=t_*\delta_{1-A}+(\gamma/q)\delta_0$. Here $q=p$, so $\mathbb EY=t_*(1-A)=q=p$ and $\Gamma(Y)=t_*(\gamma/q)(1-A)=\gamma$. Symmetry and $h(0)=0$ imply $\mathbb Eh(Y)=t_*h(1-A)=t_*h(A)$, which is again the right-hand side of~\eqref{eq:general-moment-sharp}. These computations also cover $\gamma=0$, $\gamma=p(1-p)$ and $p=0,1$, because $q-\gamma>0$ throughout and zero probabilities are allowed.
\end{proof}

To apply Lemma~\ref{lem:general-moment} to the Kleitman--West problem, we will first shift the family and approximate its normalized adjacency count by a coordinate integral. The functions needed to estimate this integral will then be verified to satisfy the hypotheses of the lemma.

\section{The dense asymptotic formula}\label{sec:dense}

A family $\mathcal F\subseteq\binom{[n]}r$ is called shifted if
\[
(F\setminus\{j\})\cup\{i\}\in\mathcal F
\]
whenever $F\in\mathcal F$, $j\in F$, $i\notin F$ and $i<j$. Thus, if $\mathcal F$ is shifted and $F\in\mathcal F$, replacing an element of $F$ by a smaller element outside $F$ gives another set in $\mathcal F$. Such families arise from the shifting method, introduced by Erd\H{o}s, Ko and Rado~\cite{EKR1961} in their study of intersecting families. Since then, shifting has been applied to a wide range of problems in extremal set theory; see, for example,~\cite{AhlswedeKhachatrian1997,HuangLohSudakov2012,FranklKupavskii2022,HuangLinialNavesPeledSudakov2016}. Here we describe the shifting operation and prove that it preserves the family size without decreasing the number of adjacent pairs. We then express the adjacency count of a shifted family as a sum of element labels, which will allow us to estimate this count by an integral.

For $i<k$, the $(i,k)$-shift, also called the $(i,k)$-compression, replaces $F\in\mathcal F$ by $(F\setminus\{k\})\cup\{i\}$ if $k\in F$, $i\notin F$ and the replacement is not already in $\mathcal F$; all other members remain unchanged. A family is shifted precisely when every such operation fixes it. This is the standard shifting operation surveyed by Frankl~\cite{Frankl1987}. To verify its effect on the number of adjacent pairs, we use the codegree identities
\begin{equation}\label{eq:codegrees}
A(\mathcal F)=\sum_{S\in\binom{[n]}{r-1}}\binom{d_{\mathcal F}(S)}2,
\qquad
\sum_{S\in\binom{[n]}{r-1}}d_{\mathcal F}(S)=r|\mathcal F|,
\end{equation}
where $d_{\mathcal F}(S)$ denotes the codegree of $S$ in $\mathcal F$, namely the number of members of $\mathcal F$ containing $S$. The first identity counts an adjacent pair by its unique intersection, and the second counts the $(r-1)$-subsets of each member of $\mathcal F$.

For a fixed family size, the following lemma shows that it suffices to consider shifted families when maximizing the number of adjacent pairs, and expresses their adjacency counts as sums over their members. This result was obtained by Harper~\cite{Harper1991} and, independently, by Ahlswede and Cai~\cite{AhlswedeCai1999}. We include a proof for completeness.

\begin{lemma}\label{lem:compression}
For every $\mathcal F\subseteq\binom{[n]}r$, there is a shifted family $\mathcal F'$ of the same size with $A(\mathcal F')\ge A(\mathcal F)$. If $\mathcal F$ is shifted, then
\begin{equation}\label{eq:rank}
A(\mathcal F)=\sum_{\{a_1<\cdots<a_r\}\in\mathcal F}\left(\sum_{\ell=1}^r a_\ell-\frac{r(r+1)}2\right).
\end{equation}
\end{lemma}

\begin{proof}
For the first assertion, start with $\mathcal F_0=\mathcal F$. Given $\mathcal F_t$, check whether there are $1\le i<k\le n$ and $F\in\mathcal F_t$ such that $k\in F$, $i\notin F$ and $(F\setminus\{k\})\cup\{i\}\notin\mathcal F_t$. If no such pair exists, stop. Otherwise, among all pairs for which such a member exists, choose the smallest $i$ and then the smallest $k$ for that $i$. This choice is possible because there are only finitely many pairs. Apply the $(i,k)$-compression to obtain $\mathcal F_{t+1}$. More explicitly, replace every member satisfying the three replacement conditions by $(F\setminus\{k\})\cup\{i\}$ and leave every other member unchanged. All conditions are checked in $\mathcal F_t$, and all replacements are made simultaneously.

Each new member lies outside $\mathcal F_t$, so it cannot coincide with a member that remains unchanged. Two distinct replaced members cannot produce the same new member either, since removing $i$ and inserting $k$ recovers the original member uniquely. Every original member therefore has a distinct resulting member, and hence $|\mathcal F_{t+1}|=|\mathcal F_t|$. Every resulting member still has size $r$, since each replacement removes one element and inserts an element that was not already present.

To compare the adjacency counts, write $d_t(S)=d_{\mathcal F_t}(S)$ for every stage $t$ and every $(r-1)$-set $S$. Two distinct members containing $S$ have intersection exactly $S$, and every adjacent pair has a unique intersection of size $r-1$. Thus the first identity in~\eqref{eq:codegrees} takes the form $A(\mathcal F_t)=\sum_{S\in\binom{[n]}{r-1}}\binom{d_t(S)}2$. We compare these summands before and after the chosen $(i,k)$-compression.

If $S$ contains neither $i$ nor $k$, a replaced member contains $S$ before the replacement if and only if its resulting member contains $S$, because the replacement changes only $i$ and $k$. Together with the correspondence between members established above, this proves $d_{t+1}(S)=d_t(S)$. If $S$ contains both $i$ and $k$, every member containing $S$ remains unchanged because it already contains $i$. No new member contains $k$, so no new member contains $S$. Hence $d_{t+1}(S)=d_t(S)$ in this case as well.

The remaining $(r-1)$-sets contain exactly one of $i,k$. Pair them as $S_i=T\cup\{i\}$ and $S_k=T\cup\{k\}$, where $T\subseteq[n]\setminus\{i,k\}$ and $|T|=r-2$. Every remaining $(r-1)$-set belongs to exactly one such pair. When $r=2$, we take $T=\varnothing$.

Fix one such $T$. Every replaced member contains $k$ but not $i$. If it contains $T$, then before replacement it contains $S_k$ but not $S_i$, whereas afterwards it contains $S_i$ but not $S_k$. This replacement therefore increases the number of members containing $S_i$ by one and decreases the number containing $S_k$ by one. If the replaced member does not contain $T$, its resulting member still does not contain $T$, because neither $i$ nor $k$ belongs to $T$. Such a replacement changes neither of the two numbers. Members that are not replaced remain unchanged. Since distinct replaced members produce distinct new members, the number of members containing $S_i$ increases by exactly the number by which the number of members containing $S_k$ decreases. Consequently, $d_{t+1}(S_i)+d_{t+1}(S_k)=d_t(S_i)+d_t(S_k)$.

We also compare the new codegree of $S_i$ with each of the two original codegrees. Every original member containing $S_i$ contains $i$ and therefore remains unchanged, so $d_{t+1}(S_i)\ge d_t(S_i)$.

To compare $d_{t+1}(S_i)$ with $d_t(S_k)$, write each original member containing $S_k$ uniquely as $T\cup\{k,x\}$, where $x\notin T\cup\{k\}$. The additional element $x$ may equal $i$, even though $S_k$ itself does not contain $i$. If $x=i$, the member is $T\cup\{i,k\}$, which remains unchanged and already contains $S_i$. If $x\ne i$, consider $T\cup\{i,x\}$. If this set already belongs to $\mathcal F_t$, it remains unchanged because it contains $i$. If it does not belong to $\mathcal F_t$, the member $T\cup\{k,x\}$ is replaced by it. In either case, $T\cup\{i,x\}$ belongs to $\mathcal F_{t+1}$ and contains $S_i$. Distinct choices of $x\ne i$ give distinct resulting sets. None of those sets contains $k$, so none coincides with the member $T\cup\{i,k\}$ obtained when $x=i$. Thus every original member containing $S_k$ determines a different member of $\mathcal F_{t+1}$ containing $S_i$, which proves $d_{t+1}(S_i)\ge d_t(S_k)$. Combining the two comparisons, we obtain $d_{t+1}(S_i)\ge\max\{d_t(S_i),d_t(S_k)\}$. Since $d_{t+1}(S_i)+d_{t+1}(S_k)=d_t(S_i)+d_t(S_k)$, we obtain
\[
\begin{aligned}
&\binom{d_{t+1}(S_i)}2+\binom{d_{t+1}(S_k)}2-\binom{d_t(S_i)}2-\binom{d_t(S_k)}2\\
&\quad=\frac12\bigl(d_{t+1}(S_i)^2+d_{t+1}(S_k)^2-d_t(S_i)^2-d_t(S_k)^2\bigr)\\
&\quad=\bigl(d_{t+1}(S_i)-d_t(S_i)\bigr)\bigl(d_{t+1}(S_i)-d_t(S_k)\bigr)\ge0.
\end{aligned}
\]
The last inequality follows because both factors are nonnegative. The summands indexed by sets containing neither or both of $i,k$ are unchanged, while the sum of the two summands indexed by $S_i,S_k$ does not decrease for any $T$. Since these cases include every $(r-1)$-set exactly once, summing proves $A(\mathcal F_{t+1})\ge A(\mathcal F_t)$.

To prove that the compression process stops, consider the sum of all element labels over all members, $W(\mathcal F_t)=\sum_{F\in\mathcal F_t}\sum_{x\in F}x$. This is a nonnegative integer. The chosen pair $i<k$ changes at least one member, and each replaced member has its element sum decreased by $k-i$. All other members remain unchanged. Therefore, $W(\mathcal F_{t+1})\le W(\mathcal F_t)-(k-i)\le W(\mathcal F_t)-1$. After $t$ steps this implies $W(\mathcal F_t)\le W(\mathcal F_0)-t$. The process cannot continue indefinitely, since a step with $t>W(\mathcal F_0)$ would make this nonnegative integer negative. It therefore stops after some finite number $s$ of steps.

Set $\mathcal F'=\mathcal F_s$. For every $F\in\mathcal F'$ and every $i<k$ with $k\in F$ and $i\notin F$, the stopping condition implies $(F\setminus\{k\})\cup\{i\}\in\mathcal F'$. Otherwise, the pair $i,k$ would still allow a compression that changes the family, contrary to termination. This is precisely the defining property of a shifted family. Since every step preserves the family size and does not decrease the adjacency count, $|\mathcal F'|=|\mathcal F|$ and $A(\mathcal F')\ge A(\mathcal F)$, proving the first assertion.

For the second assertion, suppose that $\mathcal F$ is shifted. For each $F\in\mathcal F$, write $F=\{a_1<\cdots<a_r\}$, and adopt the convention $[0]=\varnothing$. Counting adjacent pairs by replacing an element with a smaller element, we have
\[
A(\mathcal F)
=\sum_{F\in\mathcal F}\sum_{\ell=1}^r
\left|\left\{x\in[a_\ell-1]\setminus F:(F\setminus\{a_\ell\})\cup\{x\}\in\mathcal F\right\}\right|.
\]
Every $x$ counted on the right determines a member $G=(F\setminus\{a_\ell\})\cup\{x\}$ of $\mathcal F$. Since $x\notin F$, the members $F,G$ are distinct and have intersection $F\setminus\{a_\ell\}$, so they are adjacent.

Conversely, for any adjacent pair $\{F,G\}$, the sets $F\setminus G$ and $G\setminus F$ each contain a unique element. Denote these elements by $u$ and $v$, respectively, and name the two members so that $u>v$. The pair is counted in the sum for $F$ at the unique index $\ell$ with $a_\ell=u$ and with $x=v$. The removed and inserted elements are uniquely determined by $F,G$, so this pair occurs only once in the sum for $F$. It does not occur in the sum for $G$, since obtaining $F$ from $G$ requires inserting the larger element $u$ in place of $v$. Hence the displayed sum counts every adjacent pair exactly once.

Since $\mathcal F$ is shifted, $(F\setminus\{a_\ell\})\cup\{x\}\in\mathcal F$ for every $F\in\mathcal F$, every $1\le\ell\le r$ and every $x\in[a_\ell-1]\setminus F$. The membership restriction in the preceding sum is therefore automatic, and
\[
A(\mathcal F)=\sum_{F\in\mathcal F}\sum_{\ell=1}^r\bigl|[a_\ell-1]\setminus F\bigr|.
\]
The elements of $F$ smaller than $a_\ell$ are precisely its first $\ell-1$ elements, so $\bigl|[a_\ell-1]\cap F\bigr|=\ell-1$. Thus $\bigl|[a_\ell-1]\setminus F\bigr|=(a_\ell-1)-(\ell-1)$, and substitution gives
\[
\begin{aligned}
A(\mathcal F)
&=\sum_{F\in\mathcal F}\sum_{\ell=1}^r\bigl((a_\ell-1)-(\ell-1)\bigr)\\
&=\sum_{F\in\mathcal F}\left(\sum_{\ell=1}^r a_\ell-\sum_{\ell=1}^r\ell\right)\\
&=\sum_{\{a_1<\cdots<a_r\}\in\mathcal F}\left(\sum_{\ell=1}^r a_\ell-\frac{r(r+1)}2\right),
\end{aligned}
\]
proving~\eqref{eq:rank}.
\end{proof}

To relate the counting formula to an integral, we use the approach of Harper~\cite{Harper1991,Harper2004} (see also~\cite{AhlswedeCai2006}), which associates a subset of $[0,1]^r$ with a shifted family. A measurable set $D\subseteq[0,1]^r$ is called a downset if, whenever $x=(x_1,\ldots,x_r)\in D$ and $y=(y_1,\ldots,y_r)\in[0,1]^r$ satisfy $y_i\le x_i$ for every $1\le i\le r$, we have $y\in D$. We call it symmetric if it is invariant under every permutation of its coordinates. Write $|D|$ for Lebesgue measure and set $J_r(D)=\frac2r\int_D\sum_{i=1}^r x_i\,dx$. The following construction explains why symmetric downsets arise from shifted families.

\begin{lemma}\label{lem:cells}
Let $\mathcal F\subseteq\binom{[n]}r$ be shifted, let $m=|\mathcal F|$, and put $p=r!m/n^r$. There is a symmetric downset $D\subseteq[0,1]^r$ such that $|D|=p+O_r(n^{-1})$ and
\begin{equation}\label{eq:cell-moment}
J_r(D)=\frac{2(r-1)!}{n^{r+1}}A(\mathcal F)+O_r(n^{-1}).
\end{equation}
Both errors are uniform in $\mathcal F$.
\end{lemma}

\begin{proof}
If $\mathcal F=\varnothing$, take $D=\varnothing$. Otherwise, partition $(0,1]$ into the intervals $I_a=((a-1)/n,a/n]$, $1\le a\le n$, and let $\mathcal T$ consist of all ordered tuples $\mathbf a=(a_1,\ldots,a_r)$ of distinct elements of $[n]$ such that $\{a_1,\ldots,a_r\}\in\mathcal F$. Each member of $\mathcal F$ has exactly $r!$ orderings, so $|\mathcal T|=r!m$. For $\mathbf a\in\mathcal T$, write $C(\mathbf a)=I_{a_1}\times\cdots\times I_{a_r}$, and let $D_0=\bigcup_{\mathbf a\in\mathcal T}C(\mathbf a)$. These cubes are pairwise disjoint and each has volume $n^{-r}$, whence $|D_0|=r!m/n^r=p$.

To obtain a downset containing $D_0$, we define
\[
D=\bigcup_{\mathbf a\in\mathcal T}\prod_{i=1}^r[0,a_i/n].
\]
The set $D$ is measurable because it is a finite union of closed boxes. If $x\in D$ and $0\le y_i\le x_i$ for every $i$, then any box in this union containing $x$ also contains $y$, so $D$ is a downset. Moreover, $\mathcal T$ contains every coordinate permutation of each of its tuples, and hence $D$ is symmetric. Since $D_0\subseteq D$, we can split the integral defining $J_r(D)$ as
\[
J_r(D)
=\frac2r\int_D\sum_{i=1}^r x_i\,dx=\frac2r\sum_{\mathbf a\in\mathcal T}\int_{C(\mathbf a)}\sum_{i=1}^r x_i\,dx+\frac2r\int_{D\setminus D_0}\sum_{i=1}^r x_i\,dx.
\]
For the first term, integration in one coordinate leaves a factor $n^{-(r-1)}$ from the other coordinate intervals, so
\[
\int_{C(\mathbf a)}\sum_{i=1}^r x_i\,dx=\frac1{n^{r-1}}\sum_{i=1}^r\int_{(a_i-1)/n}^{a_i/n}t\,dt=\frac1{n^{r+1}}\left(\sum_{i=1}^r a_i-\frac r2\right).
\]
Each member $\{a_1<\cdots<a_r\}\in\mathcal F$ appears in $\mathcal T$ in exactly $r!$ orders, and all these orders have the same coordinate sum. Thus the sum over $\mathcal T$ is $r!$ times the sum over the increasing representations of members of $\mathcal F$. Applying~\eqref{eq:rank} and using $r(r+1)/2-r/2=r^2/2$ and $p=r!m/n^r$, we obtain
\begin{equation}\label{eq:cell-exact}
\begin{aligned}
J_r(D)
&=\frac{2(r-1)!}{n^{r+1}}\sum_{\{a_1<\cdots<a_r\}\in\mathcal F}\left(\sum_{i=1}^r a_i-\frac r2\right)+\frac2r\int_{D\setminus D_0}\sum_{i=1}^r x_i\,dx\\
&=\frac{2(r-1)!}{n^{r+1}}\left(A(\mathcal F)+\frac{r^2m}{2}\right)+\frac2r\int_{D\setminus D_0}\sum_{i=1}^r x_i\,dx\\
&=\frac{2(r-1)!}{n^{r+1}}A(\mathcal F)+\frac{rp}{n}+\frac2r\int_{D\setminus D_0}\sum_{i=1}^r x_i\,dx.
\end{aligned}
\end{equation}
To bound the volume of $D\setminus D_0$, we first identify the regions that contain its points with positive coordinates.

\begin{claim}\label{clm:cell-collisions}
We have
\[
(D\setminus D_0)\cap(0,1]^r
\subseteq
\bigcup_{1\le i<j\le r}\ \bigcup_{a=1}^n
\{y\in[0,1]^r:y_i,y_j\in I_a\}.
\]
\end{claim}
\begin{proof}[Proof of Claim~\ref{clm:cell-collisions}]
Take $y=(y_1,\ldots,y_r)\in(D\setminus D_0)\cap(0,1]^r$ and put $b_i=\lceil ny_i\rceil$, so that $y_i\in I_{b_i}$. By the construction of $D$, some $\mathbf a=(a_1,\ldots,a_r)\in\mathcal T$ satisfies $y_i\le a_i/n$, and hence $b_i\le a_i$, for every $i$. Suppose that $b_1,\ldots,b_r$ are distinct, and write $\alpha_1<\cdots<\alpha_r$ and $\beta_1<\cdots<\beta_r$ for the increasing rearrangements of $a_1,\ldots,a_r$ and $b_1,\ldots,b_r$, respectively. Since $b_i\le a_i$, every index with $a_i\le\alpha_j$ also satisfies $b_i\le\alpha_j$, and therefore
\[
\bigl|\{i\in[r]:b_i\le\alpha_j\}\bigr|
\ge\bigl|\{i\in[r]:a_i\le\alpha_j\}\bigr|
=j.
\]
Thus at least $j$ of the $b_i$ are at most $\alpha_j$, so $\beta_j\le\alpha_j$.

To show that $\{b_1,\ldots,b_r\}\in\mathcal F$, set $F_j=\{\beta_1,\ldots,\beta_j,\alpha_{j+1},\ldots,\alpha_r\}$ for $0\le j\le r$. In particular, $F_0=\{\alpha_1,\ldots,\alpha_r\}\in\mathcal F$. Suppose that $F_{j-1}\in\mathcal F$. If $\beta_j=\alpha_j$, then $F_j=F_{j-1}\in\mathcal F$. Otherwise, $\beta_j<\alpha_j$ is larger than all the previously inserted elements $\beta_1,\ldots,\beta_{j-1}$ and smaller than all the remaining elements $\alpha_j,\ldots,\alpha_r$. Hence $\beta_j\notin F_{j-1}$ and $\alpha_j\in F_{j-1}$. Since $\mathcal F$ is shifted,
\[
F_j=\bigl(F_{j-1}\setminus\{\alpha_j\}\bigr)\cup\{\beta_j\}\in\mathcal F.
\]
By induction, $\{b_1,\ldots,b_r\}=F_r\in\mathcal F$, and consequently $(b_1,\ldots,b_r)\in\mathcal T$. Since $y_i\in I_{b_i}$ for every $i$, it follows that
\[
y\in I_{b_1}\times\cdots\times I_{b_r}
=C\bigl((b_1,\ldots,b_r)\bigr)
\subseteq D_0,
\]
contrary to $y\notin D_0$. Thus there are $i<j$ and $a\in[n]$ with $b_i=b_j=a$, and hence $y_i,y_j\in I_a$. This proves the inclusion.
\end{proof}

Since the set of points in $[0,1]^r$ with at least one zero coordinate has Lebesgue measure zero, we have $|D\setminus D_0|=\bigl|(D\setminus D_0)\cap(0,1]^r\bigr|$. For fixed $i<j$ and $a$, the region where both coordinates lie in $I_a$ has volume $n^{-2}$, since the remaining coordinates range over $[0,1]$. By Claim~\ref{clm:cell-collisions} and the union bound,
\[
0\le |D|-p
=|D\setminus D_0|\le\sum_{1\le i<j\le r}\sum_{a=1}^n\bigl|\{y\in[0,1]^r:y_i,y_j\in I_a\}\bigr|
=\binom r2\,n\cdot\frac1{n^2}
=\frac{\binom r2}{n}.
\]
This proves $|D|=p+O_r(n^{-1})$. Finally, since $0\le x_i\le1$ for every $1\le i\le r$ on $[0,1]^r$, we have $0\le(2/r)\sum_{i=1}^r x_i\le(2/r)\cdot r=2$. Moreover, $p=|D_0|\le1$. Combining these bounds with~\eqref{eq:cell-exact}, we obtain
\[
\begin{aligned}
0&\le J_r(D)-\frac{2(r-1)!}{n^{r+1}}A(\mathcal F)
=\frac{rp}{n}+\frac2r\int_{D\setminus D_0}\sum_{i=1}^r x_i\,dx\\
&\le\frac r n+2|D\setminus D_0|
\le\frac r n+\frac{2\binom r2}{n}
=\frac{r^2}{n}.
\end{aligned}
\]
The two error bounds depend only on $r$ and $n$, so both are uniform in $\mathcal F$.
\end{proof}

To bound $J_r(D)$ by induction on $r$, we will use the volumes of the $(r-1)$-dimensional slices of $D$. Let $X$ denote the volume of a slice obtained by fixing the last coordinate at a uniformly chosen value in $[0,1]$, and let $X'$ be an independent copy. The slice calculation below shows that $J_r(D)=\mathbb E\min(X,X')$. Thus $J_r(D)=p-\Gamma(X)$, where $p=\mathbb EX=|D|$. This expresses $J_r(D)$ in terms of the two parameters appearing in Lemma~\ref{lem:general-moment}. To apply that lemma to slices, the relevant function is the difference between the density and the continuous extremal value. For an integer $d\ge1$ and $0\le x\le1$, define
\begin{equation}\label{eq:h-definition}
h_d(x)=x-F_d(x)
=x-\max\{x^{1+1/d},\ 2x-1+(1-x)^{1+1/d}\}.
\end{equation}
The next lemma evaluates the lower bound in Lemma~\ref{lem:general-moment} for these functions.

\begin{lemma}\label{lem:moment}
Let $d\ge1$ be an integer, let $0\le X\le1$, and put $p=\mathbb E X$, $q=\max\{p,1-p\}$ and $\gamma=\Gamma(X)$. Then
\begin{equation}\label{eq:moment-sharp}
\mathbb E h_d(X)\ge q-\frac{q^{1+2/d}}{(q-\gamma)^{1/d}}.
\end{equation}
For every admissible pair $p\in[0,1]$ and $\gamma\in[0,p(1-p)]$, equality is attained by a law supported on $\{a,1\}$ or $\{0,a\}$ for some $a\in[0,1]$, with zero masses allowed.
\end{lemma}

\begin{proof}
We verify the hypotheses of Lemma~\ref{lem:general-moment} for $h=h_d$. Put $\alpha=1+1/d\in(1,2]$. Both functions in the maximum defining $F_d$ are continuous and convex, so their maximum is continuous and convex. Hence $h_d=x-F_d$ is continuous and concave. The expression $h_d(x)=\min\{x-x^\alpha,(1-x)-(1-x)^\alpha\}$ shows that $h_d(x)=h_d(1-x)$ and $h_d(0)=h_d(1)=0$.

To identify $h_d$ on $[0,1/2]$, set $w(x)=2x-1+(1-x)^\alpha-x^\alpha$. Since $1<\alpha\le2$ and $x\le1-x$ on this interval,
\[
w''(x)=\alpha(\alpha-1)\bigl((1-x)^{\alpha-2}-x^{\alpha-2}\bigr)\le0
\qquad(0<x<1/2).
\]
Together with $w(0)=w(1/2)=0$, this proves $w(x)\ge0$ on $[0,1/2]$. Thus the second expression in the maximum defining $F_d$ is at least the first there, and
\begin{equation}\label{eq:h-branch}
h_d(x)=(1-x)-(1-x)^{1+1/d}
\qquad(0\le x\le1/2).
\end{equation}
A direct calculation using~\eqref{eq:h-branch} shows that $c_d(x)=h_d(x)/\sqrt{x(1-x)}$ is nondecreasing on $(0,1/2]$.

We may now apply Lemma~\ref{lem:general-moment}. Put $K=p(1-p)-\gamma$ and $A=K/(q-\gamma)$. Since $p(1-p)=q(1-q)$ and $0\le\gamma\le p(1-p)$, we have $q^2\le q-\gamma\le q$, and hence $1/2\le q\le1-A=q^2/(q-\gamma)\le1$. Thus $0\le A\le1/2$, and equation~\eqref{eq:h-branch} evaluates the resulting bound as
\[
\mathbb E h_d(X)\ge\frac{q-\gamma}{q}h_d(A)
=\frac{q-\gamma}{q}\bigl((1-A)-(1-A)^{1+1/d}\bigr)
=q-\frac{q^{1+2/d}}{(q-\gamma)^{1/d}}.
\]
The laws in~\eqref{eq:general-attaining-law} attain this value and have the stated support.
\end{proof}

The same verification applies to every real $d\ge1$ if $h_d$ is defined by the explicit maximum in~\eqref{eq:h-definition}. Only integer values are needed for the dimension induction.

For the dimension induction, we use Lemma~\ref{lem:moment} when $\mathbb E h_d(X)\le\Gamma(X)$. Under this additional hypothesis, the sharp estimate has the following simpler consequence.

\begin{lemma}\label{lem:scalar}
Let $d\ge1$ be an integer and let $X$ take values in $[0,1]$. If $\mathbb E h_d(X)\le\Gamma(X)$, then $\Gamma(X)\ge h_{d+1}(\mathbb E X)$.
\end{lemma}

\begin{proof}
Put $p=\mathbb EX$, $q=\max\{p,1-p\}$ and $\gamma=\Gamma(X)$. By~\eqref{eq:gini-range}, $q-\gamma\ge q^2>0$. Lemma~\ref{lem:moment} and the hypothesis imply
\[
q-\frac{q^{1+2/d}}{(q-\gamma)^{1/d}}
\le\mathbb E h_d(X)\le\gamma.
\]
Rearranging yields $\gamma\ge q-q^{1+1/(d+1)}$. By~\eqref{eq:h-definition} and the choice of $q$,
\[
h_{d+1}(p)
=\min\left\{p-p^{1+1/(d+1)},\,
(1-p)-(1-p)^{1+1/(d+1)}\right\}
\le q-q^{1+1/(d+1)}\le\gamma. \qedhere
\]
\end{proof}

We apply Lemma~\ref{lem:scalar} to the volumes of slices to prove the following bound by induction on the dimension. This bound extends the three-dimensional conclusion obtained by Ahlswede and Cai~\cite{AhlswedeCai2006}.

\begin{lemma}\label{lem:downset}
If $D\subseteq[0,1]^r$ is a measurable symmetric downset and $|D|=p$, then $J_r(D)\le F_r(p)$. Equality is attained by at least one of the symmetric downsets $D'=[0,p^{1/r}]^r$ and $D''=\{x\in[0,1]^r:\min_{1\le i\le r}x_i\le1-(1-p)^{1/r}\}$.
\end{lemma}

\begin{proof}
In dimension one, downward closure and $|D|=p$ imply $[0,p)\subseteq D\subseteq[0,p]$. Hence $J_1(D)=2\int_0^p t\,dt=p^2=F_1(p)$. Suppose $r\ge2$ and the assertion is known in dimension $r-1$. For $0\le t\le1$, let $D_t=\{(x_1,\ldots,x_{r-1}):(x_1,\ldots,x_{r-1},t)\in D\}$, and write $f(t)=|D_t|$ for its $(r-1)$-dimensional volume. These slices are symmetric downsets for almost every $t$, and $f$ is nonincreasing. Let $U,U'$ be independent uniform random variables on $[0,1]$, and put $X=f(U)$ and $X'=f(U')$. Since $U$ and $U'$ are uniformly distributed on $[0,1]$, we have $\mathbb EX=\mathbb EX'=\int_0^1f(t)\,dt$. By Fubini's theorem,
\[
\int_0^1f(t)\,dt=\int_0^1|D_t|\,dt=|D|=p.
\]
Thus $\mathbb EX=\mathbb EX'=p$. The following identities express $J_r(D)$ in terms of $X$ and $X'$.

\begin{claim}\label{clm:slice-identities}
\begin{equation}\label{eq:slice-identities}
J_r(D)=2\int_0^1tf(t)\,dt
=\mathbb E\min(X,X')
=p-\tfrac12\mathbb E|X-X'|
=p-\Gamma(X).
\end{equation}
\end{claim}

\begin{proof}[Proof of Claim~\ref{clm:slice-identities}]
By symmetry of $D$ and Fubini's theorem,
\[
J_r(D)=\frac2r\int_D\sum_{i=1}^r x_i\,dx
=2\int_Dx_r\,dx
=2\int_0^1t|D_t|\,dt
=2\int_0^1tf(t)\,dt.
\]
To calculate $\mathbb E\min(X,X')$, first note that $f$ is nonincreasing, so
\[
\min(X,X')
=\min\{f(U),f(U')\}
=f(\max\{U,U'\}).
\]
We compute the distribution of $\max\{U,U'\}$. For $0\le t\le1$, the maximum is at most $t$ precisely when both $U$ and $U'$ are at most $t$. Independence and uniformity therefore imply
\[
\mathbb P(\max\{U,U'\}\le t)=\mathbb P(U\le t,\ U'\le t)=\mathbb P(U\le t)\mathbb P(U'\le t)=t^2.
\]
Thus $\max\{U,U'\}$ has distribution function $t^2$ on $[0,1]$ and density $2t$. Integrating against this density, we obtain
\[
\mathbb E\min(X,X')
=\mathbb E f(\max\{U,U'\})
=\int_0^1f(t)\,2t\,dt
=2\int_0^1tf(t)\,dt.
\]
Finally, the identity $\min(X,X')=(X+X'-|X-X'|)/2$ and the equality $\mathbb EX=\mathbb EX'=p$ imply
\[
\mathbb E\min(X,X')
=\frac{\mathbb EX+\mathbb EX'-\mathbb E|X-X'|}{2}
=p-\frac12\mathbb E|X-X'|. \qedhere
\]
\end{proof}

To compare $\mathbb E h_{r-1}(X)$ with $\Gamma(X)$, we use the induction hypothesis on the slices.

\begin{claim}\label{clm:slice-induction}
We have
\[
J_r(D)\le p-\mathbb E h_{r-1}(X).
\]
\end{claim}

\begin{proof}[Proof of Claim~\ref{clm:slice-induction}]
By symmetry of $D$, the integrals $\int_Dx_i\,dx$ are equal for $1\le i\le r$. In particular, $\int_Dx_r\,dx=(r-1)^{-1}\int_D\sum_{i=1}^{r-1}x_i\,dx$. Starting from the definition of $J_r(D)$, we therefore obtain
\[
J_r(D)=\frac2r\int_D\sum_{i=1}^r x_i\,dx=\frac2r\left(1+\frac1{r-1}\right)\int_D\sum_{i=1}^{r-1}x_i\,dx=\frac2{r-1}\int_D\sum_{i=1}^{r-1}x_i\,dx.
\]
Fixing the last coordinate at $t$ leaves the integration region $D_t$ for the other coordinates. Fubini's theorem and the definition of $J_{r-1}$ thus imply
\[
J_r(D)=\int_0^1\left(\frac2{r-1}\int_{D_t}\sum_{i=1}^{r-1}x_i\,dx_1\cdots dx_{r-1}\right)dt=\int_0^1J_{r-1}(D_t)\,dt.
\]
For almost every $t$, the slice $D_t$ is a measurable symmetric downset in $[0,1]^{r-1}$ with volume $f(t)$. The induction hypothesis therefore implies $J_{r-1}(D_t)\le F_{r-1}(f(t))$. Integrating this inequality and using $F_{r-1}(x)=x-h_{r-1}(x)$, we obtain
\[
J_r(D)\le\int_0^1F_{r-1}(f(t))\,dt=\int_0^1f(t)\,dt-\int_0^1h_{r-1}(f(t))\,dt.
\]
The first integral equals $p$. Since $U$ is uniform on $[0,1]$ and $X=f(U)$, the second integral equals $\mathbb E h_{r-1}(X)$. Hence $J_r(D)\le p-\mathbb E h_{r-1}(X)$, as required.
\end{proof}

Comparing Claim~\ref{clm:slice-induction} with~\eqref{eq:slice-identities} shows that $\mathbb E h_{r-1}(X)\le\Gamma(X)$. By Lemma~\ref{lem:scalar}, we have $\Gamma(X)\ge h_r(\mathbb EX)=h_r(p)$. Hence $J_r(D)=p-\Gamma(X)\le p-h_r(p)=F_r(p)$.

It remains to show that equality is attained by $D'$ or $D''$. Both sets are symmetric downsets. The cube $D'$ has volume $(p^{1/r})^r=p$, and direct integration shows that
\[
J_r(D')=2p^{(r-1)/r}\int_0^{p^{1/r}}t\,dt=p^{1+1/r}.
\]
For the second construction, $[0,1]^r\setminus D''=(1-(1-p)^{1/r},1]^r$ has volume $1-p$, so $|D''|=p$. By the definition of $J_r$, symmetry and Fubini's theorem,
\[
\begin{aligned}
J_r(D'')
&=\frac2r\int_{D''}\sum_{i=1}^r x_i\,dx=2\int_{D''}x_r\,dx\\
&=2\int_{[0,1]^r}x_r\,dx-2\int_{(1-(1-p)^{1/r},1]^r}x_r\,dx\\
&=2\int_0^1t\,dt-2(1-p)^{(r-1)/r}\int_{1-(1-p)^{1/r}}^1t\,dt\\
&=2p-1+(1-p)^{1+1/r}.
\end{aligned}
\]
By the definition~\eqref{eq:F} of $F_r(p)$, we have $\max\{J_r(D'),J_r(D'')\}=F_r(p)$. Hence at least one of $D'$ and $D''$ attains equality.
\end{proof}

The continuous bound gives the required discrete estimate, including an error uniform in the number of members. We use $M_1(n,m)=\binom m2$ for singleton families, because any two distinct singletons have intersection of size zero.

\begin{lemma}\label{lem:dense}
For each fixed integer $j\ge1$, uniformly for integers $0\le m\le\binom nj$,
\begin{equation}\label{eq:dense}
M_j(n,m)=\frac{n^{j+1}}{2(j-1)!}F_j\!\left(\frac{j!m}{n^j}\right)+O_j(n^j).
\end{equation}
Let $a$ and $b$ be the largest nonnegative integers at most $n$ satisfying $\binom aj\le m$ and $\binom bj\le\binom nj-m$, respectively. Form one family by taking all $j$-subsets of $[a]$ and adding arbitrary missing members until its size is $m$. Form a second family by taking all $j$-subsets of $[n]$ not contained in $[b]$ and deleting arbitrary members until its size is $m$. At least one of these two families has at least $M_j(n,m)-O_j(n^j)$ adjacent pairs.
\end{lemma}

\begin{proof}
For $j=1$, every two distinct singletons are adjacent, so $M_1(n,m)=\binom m2=m^2/2-m/2$. Since $F_1(p)=p^2$ and $m\le n$, we obtain $M_1(n,m)=n^2F_1(m/n)/2+O(n)$, as required.

Suppose that $j\ge2$, and put $p=j!m/n^j$. By Lemma~\ref{lem:compression}, we may choose a shifted family $\mathcal F\subseteq\binom{[n]}j$ with $|\mathcal F|=m$ and $A(\mathcal F)=M_j(n,m)$. Indeed, compressing an extremal family preserves its size and does not decrease its adjacency count, so the resulting family is still extremal. By Lemma~\ref{lem:cells}, there is a measurable symmetric downset $D\subseteq[0,1]^j$. More precisely, the estimates in its proof show that $0\le |D|-p\le\binom j2/n$ and
\[
0\le J_j(D)-\frac{2(j-1)!}{n^{j+1}}M_j(n,m)\le\frac{j^2}{n},
\]
where we used $A(\mathcal F)=M_j(n,m)$. Applying Lemma~\ref{lem:downset}, we therefore obtain
\[
\frac{2(j-1)!}{n^{j+1}}M_j(n,m)\le J_j(D)\le F_j(|D|).
\]
To replace $F_j(|D|)$ by an expression in $p$, we compare the two expressions defining $F_j$. Since $0\le p\le |D|\le1$ and $1+1/j\le2$,
\[
|D|^{1+1/j}-p^{1+1/j}=\left(1+\frac1j\right)\int_p^{|D|}t^{1/j}\,dt\le2(|D|-p).
\]
For the other expression, $1-|D|\le1-p$ implies
\[
\begin{aligned}
&\bigl(2|D|-1+(1-|D|)^{1+1/j}\bigr)-\bigl(2p-1+(1-p)^{1+1/j}\bigr)\\
&\qquad=2(|D|-p)+(1-|D|)^{1+1/j}-(1-p)^{1+1/j}\le2(|D|-p).
\end{aligned}
\]
Each expression defining $F_j(|D|)$ is therefore at most $F_j(p)+2(|D|-p)$. Taking their maximum and using the bound on $|D|-p$, we obtain
\[
\frac{2(j-1)!}{n^{j+1}}M_j(n,m)\le F_j(|D|)\le F_j(p)+2(|D|-p)\le F_j(p)+\frac{2\binom j2}{n}.
\]
Multiplying by $n^{j+1}/(2(j-1)!)$ now yields
\[
M_j(n,m)\le\frac{n^{j+1}}{2(j-1)!}F_j\!\left(\frac{j!m}{n^j}\right)+\frac{\binom j2}{(j-1)!}\,n^j.
\]
The coefficient of the error term depends only on $j$, which proves the required upper bound uniformly in $m$.

It remains to prove the lower bound. We give one construction for each of the two expressions defining $F_j(p)$, adjusting the number of members to be exactly $m$. For the first construction, take all $j$-subsets of $[a]$, where $a\le n$ is largest subject to $\binom aj\le m$. Since $m-\binom aj=O_j(n^{j-1})$, we have $(a/n)^j=p+O_j(n^{-1})$. Its adjacency count satisfies
\[
\begin{aligned}
\frac{j(a-j)}2\binom aj
&=\frac{j(a-j)}2\left(\frac{a^j}{j!}+O_j(n^{j-1})\right)=\frac{a^{j+1}}{2(j-1)!}+O_j(n^j)\\
&=\frac{n^{j+1}}{2(j-1)!}\left(\left(\frac an\right)^j\right)^{1+1/j}+O_j(n^j)=\frac{n^{j+1}}{2(j-1)!}p^{1+1/j}+O_j(n^j).
\end{aligned}
\]
The last equality uses the power-function estimate established above. Adding missing $j$-sets until the family has size $m$ cannot decrease its count.

For the second construction, write $N=\binom nj$, let $q=N-m$, and take $b\le n$ largest subject to $\binom bj\le q$. Let $\mathcal H=\binom{[n]}j\setminus\binom{[b]}j$ be the family before deleting surplus members. Every adjacent pair in $\mathcal H$ has a unique intersection $S$ of size $j-1$, and its two members can be written as $S\cup\{x\}$ and $S\cup\{y\}$, where $x,y\in[n]\setminus S$ are distinct.

If $S\not\subseteq[b]$, both members belong to $\mathcal H$ for every choice of $\{x,y\}$. There are $\binom n{j-1}-\binom b{j-1}$ choices for $S$, and each permits $\binom{n-j+1}{2}$ choices for $\{x,y\}$.

If $S\subseteq[b]$, both members belong to $\mathcal H$ precisely when $x,y\notin[b]$. There are $\binom b{j-1}$ choices for $S$, and each permits $\binom{n-b}{2}$ choices for $\{x,y\}$. Since each adjacent pair is counted exactly once, summing the two contributions gives
\[
A(\mathcal H)=\left(\binom n{j-1}-\binom b{j-1}\right)\binom{n-j+1}{2}+\binom b{j-1}\binom{n-b}{2}.
\]
Using $\binom n{j-1}(n-j+1)=jN$ and $\binom b{j-1}(b-j+1)=j\binom bj$, we obtain
\[
\begin{aligned}
A(\mathcal H)&=\frac{j(n-j)}2N-\frac{j(2n-b-j)}2\binom bj\\
&=\frac{j(n-j)}2N-j(n-j)\binom bj+\frac{j(b-j)}2\binom bj.
\end{aligned}
\]
The family $\mathcal H$ has $N-\binom bj=m+O_j(n^{j-1})$ members. Deleting one member removes at most $j(n-j)$ adjacent pairs, since an adjacent member is obtained by replacing one of its $j$ elements with one of the $n-j$ elements outside it. Removing the surplus therefore changes the count by at most $j(n-j)O_j(n^{j-1})=O_j(n^j)$. Since $j!q/n^j=1-p+O_j(n^{-1})$, expansion of the preceding expression shows that the adjusted family has
\[
\frac{n^{j+1}}{2(j-1)!}\left(2p-1+(1-p)^{1+1/j}\right)+O_j(n^j)
\]
adjacent pairs. Taking the better of the two constructions proves the lower bound and, together with the upper bound, the assertion about their adjacency counts. When $m=0$, the empty family has no adjacent pairs. When $m=\binom nj$, the complete family has $j(n-j)\binom nj/2$ adjacent pairs. These values also satisfy~\eqref{eq:dense}.
\end{proof}

The additive error in~\eqref{eq:dense} is negligible when $m$ is a positive fraction of $n^j$, but it does not settle all smaller sizes. At a smaller critical scale, we instead apply the dense formula after removing an appropriate common subset from each of finitely many parts of the family.

\section{Critical scales and the proof of Theorem 1.1}\label{sec:scales}

By Lemma~\ref{lem:dense}, we can estimate $M_r(n,m)$ with an additive error $O_r(n^r)$ uniformly in $m$, but this error need not be negligible when $m\asymp n^j$ with $1\le j<r$. Indeed, each member of an $r$-uniform family has at most $r(n-r)$ adjacent members, obtained by replacing one of its elements with an element outside it. Thus $A(\mathcal F)\le r(n-r)m/2=O(n^{j+1})$, whereas $n^r\ge n^{j+1}$. We therefore cannot deduce a vanishing relative error from the estimate with error $O_r(n^r)$ alone.

To apply Lemma~\ref{lem:dense} for $j$-element sets, we first consider subfamilies whose members contain a common $(r-j)$-set. We will verify that deleting this subset preserves both the number of members and the number of adjacent pairs within each subfamily; the choice of suitable subfamilies and the bounds for the remaining contributions will be established below. More precisely, suppose that every member of a subfamily $\mathcal F_i\subseteq\binom{[n]}r$ contains the same $(r-j)$-set $T_i$, and put $\mathcal H_i=\{E\setminus T_i:E\in\mathcal F_i\}$. Each member of $\mathcal H_i$ has $r-(r-j)=j$ elements. Since adjoining $T_i$ recovers the original member, the deletion is one-to-one, so $|\mathcal H_i|=|\mathcal F_i|$. For distinct $E,E'\in\mathcal F_i$, we also have
\[
|(E\setminus T_i)\cap(E'\setminus T_i)|=|E\cap E'|-(r-j).
\]
Consequently, $|E\cap E'|=r-1$ if and only if the two sets obtained by deletion intersect in $j-1$ elements. The deletion therefore preserves every adjacent pair, and $A(\mathcal H_i)=A(\mathcal F_i)$. Writing $m_i=|\mathcal F_i|$ and viewing $\mathcal H_i$ as a $j$-uniform family on $[n]$, we apply Lemma~\ref{lem:dense} to obtain
\[
A(\mathcal F_i)=A(\mathcal H_i)\le M_j(n,m_i)=\frac{n^{j+1}}{2(j-1)!}F_j\!\left(\frac{j!m_i}{n^j}\right)+O_j(n^j).
\]
The error term is uniform in $m_i$, so the estimate remains valid even when $m_i=o(n^j)$. For a bounded number of such parts, the sum of the errors is still $O(n^j)$ and is therefore negligible compared with $n^{j+1}$.

We will show below how to choose disjoint subfamilies $\mathcal F_1,\ldots,\mathcal F_K$ with these fixed subsets, leaving the remaining members in $\mathcal F_0$. For a fixed upper bound on $m/n^j$ and a fixed error tolerance, $K$ will be bounded independently of $n$. We will bound the number of adjacent pairs inside $\mathcal F_0$ and between different parts by $(j-1)nm_0/2$ plus an error that can be made arbitrarily small compared with $n^{j+1}$, where $m_0=|\mathcal F_0|$. Before constructing the parts and proving these estimates, we explain how to combine their contributions into a bound depending only on $m$.

Write $\beta_j=(j-1)/j$ for each integer $j\ge1$, and put $y_i=j!m_i/n^j$ for $0\le i\le K$. For $i\ge1$, the family obtained by deletion has at most $\binom nj$ members, so $0\le y_i\le1$; for the remaining part, we only require $y_0\ge0$. Since $(j-1)nm_0/2=n^{j+1}\beta_jy_0/(2(j-1)!)$, the terms to be combined have the form
\[
\frac{n^{j+1}}{2(j-1)!}\left(\beta_jy_0+\sum_{i=1}^K F_j(y_i)\right).
\]
Since $\sum_{i=0}^K y_i=j!m/n^j$, we seek an upper bound for $\beta_jy_0+\sum_{i=1}^K F_j(y_i)$ that depends only on $j!m/n^j$, rather than on the individual subfamily sizes.

Both functions whose maximum defines $F_j$ in~\eqref{eq:F} are convex, and hence $F_j$ is convex. Since $F_j$ is convex and satisfies $F_j(0)=0$ and $F_j(1)=1$, we have $F_j(u)\le(1-u)F_j(0)+uF_j(1)=u$ for $0\le u\le1$. For the lower bound, we apply the tangent-line inequality to the convex function $z\mapsto z^{1+1/j}$ at $z=1$ and obtain
\begin{equation}\label{eq:all-profile-bounds}
F_j(u)\ge 2u-1+(1-u)^{1+1/j}\ge 2u-1+1-\left(1+\frac1j\right)u=\beta_ju
\qquad(0\le u\le1).
\end{equation}
The required subfamilies and the bounds for the remaining contributions will be established in Lemma~\ref{lem:all-extraction} for $2\le j<r$, with $j=1$ treated separately in the proof of Lemma~\ref{lem:all-critical}. Before proving these assertions, we establish the inequality needed to combine the resulting bounds.

\begin{lemma}\label{lem:all-packing}
Suppose that $F:[0,1]\to\mathbb R$ is convex, $F(0)=0$, $F(1)=1$, and $F(u)\ge bu$. If $y_0\ge0$, $y_1,\ldots,y_K\in[0,1]$, and $t=\sum_{i=0}^K y_i$, then
\[
by_0+\sum_{i=1}^K F(y_i)\le\lfloor t\rfloor+F(\{t\}).
\]
\end{lemma}

\begin{proof}
To handle the term $by_0$ when $y_0$ may exceed one, write $q=\lfloor y_0\rfloor$ and $\theta=y_0-q$, so that $q$ is a nonnegative integer and $0\le\theta<1$. Since $F(u)\ge bu$ for every $u\in[0,1]$, we have $b\le F(1)$ and $b\theta\le F(\theta)$. As $q\ge0$, we obtain
\[
by_0+\sum_{i=1}^K F(y_i)=qb+b\theta+\sum_{i=1}^K F(y_i)\le qF(1)+F(\theta)+\sum_{i=1}^K F(y_i).
\]
The right-hand side is the sum of the $F$-values at $q$ copies of $1$, at $\theta$, and at $y_1,\ldots,y_K$. These numbers all lie in $[0,1]$, and their sum is $q+\theta+\sum_{i=1}^K y_i=t$.

Suppose that two numbers $u,v$ in this finite list lie in $(0,1)$, and put $s=u+v$. If $s\le1$, we use the convexity of $F$ on $[0,s]$ to obtain
\[
F(u)+F(v)\le\frac{u+v}{s}F(s)+\frac{u+v}{s}F(0)=F(s)+F(0).
\]
We may therefore replace $u,v$ by $s,0$, preserving their sum without decreasing the sum of their $F$-values. If instead $s>1$, then $s<2$ and $s-1<u,v<1$. Since $F$ is convex on $[s-1,1]$, we have $F(u)\le\frac{1-u}{2-s}F(s-1)+\frac{1-v}{2-s}F(1)$ and $F(v)\le\frac{1-v}{2-s}F(s-1)+\frac{1-u}{2-s}F(1)$. Since $(1-u)+(1-v)=2-s$, we obtain $F(u)+F(v)\le F(s-1)+F(1)$ by adding these two inequalities. Thus we may replace $u,v$ by $s-1,1$, again preserving their sum without decreasing the sum of their $F$-values.

Repeat the appropriate replacement whenever two numbers lie in $(0,1)$, leaving all other numbers unchanged. Each replacement keeps every number in $[0,1]$, preserves their total sum $t$, and decreases the number of entries in $(0,1)$ by at least one. The process therefore stops after finitely many steps with at most one number in $(0,1)$. Let $\ell$ be the number of entries equal to $1$, and let $z$ be the possible remaining entry in $(0,1)$, setting $z=0$ if no such entry remains. Then $t=\ell+z$ with $0\le z<1$, so $\ell=\lfloor t\rfloor$ and $z=\{t\}$. Since $F(0)=0$ and $F(1)=1$, we conclude from the initial estimate and the successive replacements that
\[
by_0+\sum_{i=1}^K F(y_i)\le\ell F(1)+F(z)=\ell+F(z)=\lfloor t\rfloor+F(\{t\}).
\]
This also covers integer $t$, for which $z=0$ and $F(z)=0$.
\end{proof}

We now apply Lemma~\ref{lem:all-packing} to $F_j$ and express the resulting bound in terms of $G_j$.

\begin{corollary}\label{cor:all-packing}
Fix an integer $j\ge1$, and recall that $\beta_j=(j-1)/j$. If $K\ge0$ is an integer, $y_0\ge0$, and $y_1,\ldots,y_K\in[0,1]$, then
\[
\beta_jy_0+\sum_{i=1}^K F_j(y_i)\le G_j\!\left(\sum_{i=0}^K y_i\right).
\]
\end{corollary}

\begin{proof}
As verified above, $F_j$ is convex, $F_j(0)=0$, $F_j(1)=1$, and $F_j(u)\ge\beta_ju$ for $0\le u\le1$. Applying Lemma~\ref{lem:all-packing} with $F=F_j$ and $b=\beta_j$, and using the definition~\eqref{eq:G} of $G_j$, we obtain the stated inequality.
\end{proof}

We will use Corollary~\ref{cor:all-packing} to combine the bounds for the extracted parts. To justify that only finitely many parts are needed, we require a count of configurations contained in a family of high-codegree subsets. We will obtain the required upper bound by the entropy method. We first recall the definitions and basic properties of entropy needed below. For a random variable $Y$ taking values in a finite set $\mathcal A$, its entropy is defined by
\[
\mathsf H(Y)=-\sum_{y\in\mathcal A}\mathbb P(Y=y)\log\mathbb P(Y=y),
\]
where $0\log0=0$ and logarithms are natural. If $Z$ is another random variable taking values in a finite set $\mathcal B$, the conditional entropy of $Y$ given $Z$ is
\[
\mathsf H(Y\mid Z)=-\sum_{\substack{z\in\mathcal B\\ \mathbb P(Z=z)>0}}\mathbb P(Z=z)\sum_{y\in\mathcal A}\mathbb P(Y=y\mid Z=z)\log\mathbb P(Y=y\mid Z=z).
\]
Thus conditional entropy is the average of the entropies of the conditional distributions.

We will use the following standard properties of entropy; see Cover and Thomas~\cite{CoverThomas2006}.

\begin{proposition}\label{prop:entropy-properties}
Let $Y$, $Z$ and $W$ be jointly distributed random variables taking values in finite sets, and suppose that $Y$ takes values in $\mathcal A$. Then the following statements hold.
\begin{enumerate}[label=\textup{(\roman*)}]
\item We have $0\le\mathsf H(Y)\le\log|\mathcal A|$, with equality in the upper bound if and only if $Y$ is uniformly distributed on $\mathcal A$.
\item The chain rule states that $\mathsf H(Y,Z)=\mathsf H(Z)+\mathsf H(Y\mid Z)$, where $\mathsf H(Y,Z)$ denotes the entropy of the random vector $(Y,Z)$.
\item Conditional entropy is nonnegative, so $\mathsf H(Y\mid Z)\ge0$.
\item We have $\mathsf H(Y\mid Z,W)\le\mathsf H(Y\mid Z)$, so conditioning on additional variables cannot increase conditional entropy.
\end{enumerate}
\end{proposition}

Entropy methods have been widely used to establish upper bounds in extremal and enumerative combinatorics; see, for example,~\cite{ChungGrahamFranklShearer1986,Friedgut2004,MadimanTetali2010,BalisterBollobas2012,AtseriasGroheMarx2013,Simkin2023,ChaoYu2026}. In particular, Friedgut and Kahn~\cite{FriedgutKahn1998} used fractional edge covers and Shearer's entropy inequality to bound the number of copies of one hypergraph in another. The estimate needed here is the corresponding bound for injective embeddings, with an explicit constant accounting for the possible orderings of the vertices in an edge.

To prove this estimate, we choose an embedding uniformly at random and use the fractional covering assumptions to bound its entropy by the weighted entropies of its edge images.

\begin{lemma}\label{lem:all-entropy}
Let $s\ge1$ and let $H$ be a fixed nonempty $s$-uniform hypergraph without isolated vertices. Suppose that nonnegative numbers $\lambda_e$, indexed by $e\in E(H)$, satisfy $\sum_{e\ni v}\lambda_e\ge1$ for every $v\in V(H)$. For an $s$-uniform family $\mathcal Q$, let $\operatorname{Emb}(H,\mathcal Q)$ denote the number of injective maps from $V(H)$ to the ground set of $\mathcal Q$ that send every edge of $H$ to a member of $\mathcal Q$. Then
\[
\operatorname{Emb}(H,\mathcal Q)\le(s!|\mathcal Q|)^{\sum_{e\in E(H)}\lambda_e}.
\]
\end{lemma}

\begin{proof}
If there are no embeddings, the bound holds. Otherwise choose an embedding uniformly at random and let $X_v$ be the image of vertex $v$. Fix an order on $V(H)$. By the chain rule in Proposition~\ref{prop:entropy-properties}(ii) and the conditioning inequality in Proposition~\ref{prop:entropy-properties}(iv), we have, for every edge $e$,
\[
\mathsf H((X_v)_{v\in e})
=\sum_{v\in e}\mathsf H(X_v\mid X_u:u\in e, u<v)
\ge\sum_{v\in e}\mathsf H(X_v\mid X_u:u<v).
\]
Multiply by $\lambda_e\ge0$ and sum over $e$. Since conditional entropy is nonnegative by Proposition~\ref{prop:entropy-properties}(iii), we can apply the fractional covering assumptions and the chain rule in Proposition~\ref{prop:entropy-properties}(ii) to obtain
\[
\log\operatorname{Emb}(H,\mathcal Q)
=\mathsf H((X_v)_{v\in V(H)})
\le\sum_{e\in E(H)}\lambda_e\mathsf H((X_v)_{v\in e})
\le\sum_{e\in E(H)}\lambda_e\log(s!|\mathcal Q|).
\]
By Proposition~\ref{prop:entropy-properties}(i), the first equality holds because the embedding is chosen uniformly, and the last inequality holds because the ordered image of any edge has at most $s!|\mathcal Q|$ possibilities. Exponentiating both sides, we obtain the stated bound.
\end{proof}

To apply Lemma~\ref{lem:dense} to $j$-uniform families, we still need to extract boundedly many subfamilies whose members contain a common $(r-j)$-set and estimate the adjacent pairs not contained in any one of these subfamilies. The remaining family need not be small. Instead, we distinguish its members containing at least $j$ subsets of size $r-1$ whose codegrees in $\mathcal F$ are at least a chosen threshold. We will bound the number of these members by applying Lemma~\ref{lem:all-entropy}, while every other remaining member contains at most $j-1$ such subsets. We can then estimate the adjacent pairs in the remaining family by counting their common $(r-1)$-subsets. The following lemma establishes the required partition and bounds both the adjacent pairs within the remaining family and those between different subfamilies.

\begin{lemma}\label{lem:all-extraction}
Fix integers $r>j\ge2$, a constant $C>0$, and parameters $0<\eta,\delta<1$. Every family $\mathcal F\subseteq\binom{[n]}r$ of size $m\le Cn^j$ has a partition
\[
\mathcal F=\mathcal F_1\mathbin{\dot\cup}\cdots\mathbin{\dot\cup}\mathcal F_K\mathbin{\dot\cup}\mathcal F_0
\]
such that all members of $\mathcal F_i$ contain one fixed $(r-j)$-set $T_i$ for $i\ge1$, and $K\le r\binom{r-1}{j-1}C/(\eta\delta)$. Writing $m_i=|\mathcal F_i|$, we have
\begin{equation}\label{eq:all-extraction}
A(\mathcal F)\le\sum_{i=1}^K M_j(n,m_i)+\frac{j-1}{2}nm_0
+O_{r,j}\!\left(mn\left(\eta+\frac{\delta^{1/(j-1)}}\eta\right)+Km\right).
\end{equation}
\end{lemma}

\begin{proof}
Let $\mathcal S=\{S\in\binom{[n]}{r-1}:d_{\mathcal F}(S)\ge\eta n\}$. By \eqref{eq:codegrees}, $|\mathcal S|\le rm/(\eta n)$. For $T\in\binom{[n]}{r-j}$, write $b(T)=|\{S\in\mathcal S:T\subset S\}|$. Counting each $S\in\mathcal S$ once for each of its $(r-j)$-subsets, we obtain
\begin{equation}\label{eq:all-core-incidences}
\sum_{T\in\binom{[n]}{r-j}}b(T)
=\binom{r-1}{j-1}|\mathcal S|
\le\binom{r-1}{j-1}\frac{rm}{\eta n}.
\end{equation}
List every $T$ satisfying $b(T)\ge\delta n^{j-1}$ as $T_1,\ldots,T_K$. Dividing the bound in~\eqref{eq:all-core-incidences} by $\delta n^{j-1}$ and using $m\le Cn^j$, we obtain the stated bound on $K$. With this ordering fixed, for $1\le i\le K$ define
\[
\mathcal F_i=\left\{E\in\mathcal F:T_i\subseteq E,\ T_h\nsubseteq E\text{ for every }1\le h<i\right\}.
\]
The members containing none of $T_1,\ldots,T_K$ form the remaining subfamily
\[
\mathcal F_0=\left\{E\in\mathcal F:T_i\nsubseteq E\text{ for every }1\le i\le K\right\}.
\]
If $E\in\mathcal F$ contains at least one of $T_1,\ldots,T_K$, then it belongs to $\mathcal F_i$ precisely when $i$ is the smallest index for which $T_i\subseteq E$; otherwise it belongs to $\mathcal F_0$. Thus every member belongs to exactly one of these subfamilies, and $\mathcal F=\mathcal F_1\mathbin{\dot\cup}\cdots\mathbin{\dot\cup}\mathcal F_K\mathbin{\dot\cup}\mathcal F_0$.

To estimate the adjacent pairs within $\mathcal F_i$ for $i\ge1$, put $\mathcal H_i=\{E\setminus T_i:E\in\mathcal F_i\}$. Since every member of $\mathcal F_i$ contains $T_i$ and $|T_i|=r-j$, we have $\mathcal H_i\subseteq\binom{[n]\setminus T_i}{j}$. Adjoining $T_i$ recovers the original member, so $|\mathcal H_i|=|\mathcal F_i|=m_i$. For distinct $E,E'\in\mathcal F_i$, we have $T_i\subseteq E\cap E'$ and hence
\[
|(E\setminus T_i)\cap(E'\setminus T_i)|=|E\cap E'|-(r-j).
\]
Consequently, $|E\cap E'|=r-1$ if and only if $|(E\setminus T_i)\cap(E'\setminus T_i)|=j-1$, so deleting $T_i$ preserves adjacency. Since $\mathcal H_i$ can also be regarded as a $j$-uniform family on $[n]$, we obtain $A(\mathcal F_i)=A(\mathcal H_i)\le M_j(n,m_i)$.

We next bound the number of adjacent pairs whose members lie in different subfamilies. Fix $1\le i\le K$ and consider an adjacent pair with one member $E'\in\mathcal F_i$ and the other member $E\in\mathcal F_0\cup\bigcup_{\ell=i+1}^K\mathcal F_\ell$. By the definitions of these subfamilies, we have $T_i\subseteq E'$ and $T_i\nsubseteq E$, so $T_i\setminus E\subseteq E'\setminus E$. Since $E$ and $E'$ are adjacent,
\[
1\le |T_i\setminus E|\le |E'\setminus E|=r-|E\cap E'|=1.
\]
Thus $T_i\setminus E=E'\setminus E=\{x\}$, where $x$ is determined by $i$ and $E$. Every member of $\mathcal F_i$ adjacent to $E$ must therefore have the form $(E\setminus\{y\})\cup\{x\}$ for some $y\in E$, so
\[
\bigl|\{E'\in\mathcal F_i:|E\cap E'|=r-1\}\bigr|\le\bigl|\{(E\setminus\{y\})\cup\{x\}:y\in E\}\bigr|\le |E|=r.
\]
The same bound holds if $E$ has no adjacent member in $\mathcal F_i$. Since each adjacent pair between different subfamilies is counted exactly once, we bound their total number by
\[
\sum_{i=1}^K\sum_{E\in\mathcal F_0\cup\bigcup_{\ell=i+1}^K\mathcal F_\ell}\bigl|\{E'\in\mathcal F_i:|E\cap E'|=r-1\}\bigr|\le r\sum_{i=1}^K\left(m_0+\sum_{\ell=i+1}^K m_\ell\right)\le rKm.
\]

It remains to bound $A(\mathcal F_0)$. We first estimate the size of the family
\[
\mathcal D=\left\{E\in\mathcal F_0:\left|\binom E{r-1}\cap\mathcal S\right|\ge j\right\}.
\]
We call the members of $\mathcal D$ exceptional and write $D=|\mathcal D|$.

Fix an $(r-j)$-set $T\subseteq[n]$ with $T\notin\{T_1,\ldots,T_K\}$. From each set $S\in\mathcal S$ containing $T$, delete $T$, and let
\[
\mathcal S_T=\{S\setminus T:S\in\mathcal S,\ T\subseteq S\}.
\]
Each member of $\mathcal S_T$ has size $(r-1)-(r-j)=j-1$, so we regard $\mathcal S_T$ as the edge set of a $(j-1)$-uniform hypergraph on $[n]\setminus T$. Since $S=T\cup(S\setminus T)$ whenever $T\subseteq S$, the original set $S$ is uniquely determined by $S\setminus T$. Hence $|\mathcal S_T|=b(T)$. For this fixed $(r-j)$-set $T$, define $\mathcal D_T$ to consist of the exceptional members $E\in\mathcal D$ satisfying $T\subseteq E$ and $\binom{E\setminus T}{j-1}\subseteq\mathcal S_T$, so
\[
\mathcal D_T=\left\{E\in\mathcal D:T\subseteq E,\ \binom{E\setminus T}{j-1}\subseteq\mathcal S_T\right\}.
\]

We show that every member of $\mathcal D$ belongs to at least one of these families $\mathcal D_T$. Take $E\in\mathcal D$ and choose distinct sets $S_1,\ldots,S_j\in\binom E{r-1}\cap\mathcal S$. Since $|E|=r$ and $|S_\ell|=r-1$, there are pairwise distinct elements $x_1,\ldots,x_j\in E$ such that $S_\ell=E\setminus\{x_\ell\}$. Set
\[
T:=\bigcap_{\ell=1}^jS_\ell=\bigcap_{\ell=1}^j\bigl(E\setminus\{x_\ell\}\bigr)=E\setminus\{x_1,\ldots,x_j\}.
\]
Thus $|T|=r-j$ and $T\subseteq E$. Since $E\in\mathcal F_0$ contains none of $T_1,\ldots,T_K$, we have $T\notin\{T_1,\ldots,T_K\}$. Moreover, $E\setminus T=\{x_1,\ldots,x_j\}$, so
\[
\binom{E\setminus T}{j-1}=\bigl\{(E\setminus T)\setminus\{x_\ell\}:1\le\ell\le j\bigr\}=\{S_\ell\setminus T:1\le\ell\le j\}\subseteq\mathcal S_T.
\]
Therefore $E\in\mathcal D_T$. Since $E$ was arbitrary, we have $\mathcal D\subseteq\bigcup_{T\in\binom{[n]}{r-j}\setminus\{T_1,\ldots,T_K\}}\mathcal D_T$. The families $\mathcal D_T$ need not be disjoint, but the size of their union is at most the sum of their sizes. It therefore suffices to bound $|\mathcal D_T|$ for each fixed $T$.

Fix $T\in\binom{[n]}{r-j}$ with $T\notin\{T_1,\ldots,T_K\}$, and let $H_j$ be the hypergraph with vertex set $[j]$ and edge set $\binom{[j]}{j-1}$. For each $E\in\mathcal D_T$, choose a bijection $\varphi_E:[j]\to E\setminus T$. Since every edge $e$ of $H_j$ has $j-1$ vertices, we have $\varphi_E(e)\in\binom{E\setminus T}{j-1}\subseteq\mathcal S_T$. Thus $\varphi_E$ is an injective embedding of $H_j$ into $\mathcal S_T$. For fixed $T$, different members $E\in\mathcal D_T$ have different image sets $\varphi_E([j])=E\setminus T$, since $E=T\cup\varphi_E([j])$. Hence the embeddings chosen for different members are distinct. To apply Lemma~\ref{lem:all-entropy}, assign weight $\lambda_e=1/(j-1)$ to every edge of $H_j$. Each vertex belongs to exactly $j-1$ edges, so $\sum_{e\ni v}\lambda_e=1$ for every $v\in[j]$. Since $H_j$ has $j$ edges, we also have $\sum_{e\in E(H_j)}\lambda_e=j/(j-1)$. By Lemma~\ref{lem:all-entropy} with $s=j-1$, we conclude that
\[
|\mathcal D_T|=\bigl|\{\varphi_E:E\in\mathcal D_T\}\bigr|\le\operatorname{Emb}(H_j,\mathcal S_T)\le\bigl((j-1)!\,|\mathcal S_T|\bigr)^{j/(j-1)}=\bigl((j-1)!\,b(T)\bigr)^{j/(j-1)}.
\]
For $j=2$, the two edges of $H_j$ are singletons and both weights are one, so this application of Lemma~\ref{lem:all-entropy} remains valid.

By the definition of $T_1,\ldots,T_K$, we have $b(T)<\delta n^{j-1}$ for every $T$ in the union above. Consequently, $b(T)^{j/(j-1)}=b(T)b(T)^{1/(j-1)}\le\delta^{1/(j-1)}n\,b(T)$. Summing the bounds for $|\mathcal D_T|$ and applying~\eqref{eq:all-core-incidences}, we obtain
\[
\begin{aligned}
D=|\mathcal D|
&\le\sum_{T\in\binom{[n]}{r-j}\setminus\{T_1,\ldots,T_K\}}|\mathcal D_T|\le((j-1)!)^{j/(j-1)}\sum_{T\in\binom{[n]}{r-j}\setminus\{T_1,\ldots,T_K\}}b(T)^{j/(j-1)}\\
&\le((j-1)!)^{j/(j-1)}\delta^{1/(j-1)}n\sum_{T\in\binom{[n]}{r-j}}b(T)
=O_{r,j}\!\left(\frac{\delta^{1/(j-1)}}\eta m\right).
\end{aligned}
\]

Since each nonexceptional member of $\mathcal F_0$ contains at most $j-1$ sets from $\mathcal S$ and each exceptional member contains at most $r$, we obtain $\sum_{S\in\mathcal S}d_{\mathcal F_0}(S)\le(j-1)m_0+rD$ by counting pairs $(E,S)$ with $E\in\mathcal F_0$, $S\in\mathcal S$ and $S\subseteq E$. For each $S\in\binom{[n]}{r-1}$, there are $\binom{d_{\mathcal F_0}(S)}2$ pairs of distinct members of $\mathcal F_0$ containing $S$, and the intersection of each such pair is exactly $S$. Since every adjacent pair has a unique intersection of size $r-1$, we count each pair exactly once by summing over all $S\in\binom{[n]}{r-1}$. Thus, by the first identity in~\eqref{eq:codegrees},
\[
2A(\mathcal F_0)=\sum_{S\in\binom{[n]}{r-1}}d_{\mathcal F_0}(S)\bigl(d_{\mathcal F_0}(S)-1\bigr)\le\sum_{S\in\binom{[n]}{r-1}}d_{\mathcal F_0}(S)^2.
\]
For $S\in\binom{[n]}{r-1}\setminus\mathcal S$, we have $d_{\mathcal F_0}(S)\le d_{\mathcal F}(S)<\eta n$ because $\mathcal F_0\subseteq\mathcal F$ and $S\notin\mathcal S$, while $d_{\mathcal F_0}(S)\le n$ for every $S\in\binom{[n]}{r-1}$. Splitting the sum according to whether $S$ belongs to $\mathcal S$, and using $\sum_{S\in\binom{[n]}{r-1}}d_{\mathcal F_0}(S)=rm_0$, we obtain
\[
\begin{aligned}
\sum_{S\in\binom{[n]}{r-1}}d_{\mathcal F_0}(S)^2
&=\sum_{S\in\mathcal S}d_{\mathcal F_0}(S)^2+\sum_{S\in\binom{[n]}{r-1}\setminus\mathcal S}d_{\mathcal F_0}(S)^2\\
&\le n\sum_{S\in\mathcal S}d_{\mathcal F_0}(S)+\eta n\sum_{S\in\binom{[n]}{r-1}\setminus\mathcal S}d_{\mathcal F_0}(S)\\
&\le n\bigl((j-1)m_0+rD\bigr)+r\eta nm_0.
\end{aligned}
\]
Combining this bound with the estimates for the extracted parts and for cross-part pairs, we obtain~\eqref{eq:all-extraction}, as required.
\end{proof}

With $\eta$ and $\delta$ fixed, the partition has boundedly many parts, so the additive errors from Lemma~\ref{lem:dense} are negligible on the scale $n^{j+1}$. We then apply Corollary~\ref{cor:all-packing} to combine the bounds for the parts and the remainder. The linear scale requires a separate extraction because no exponent $1/(j-1)$ is available when $j=1$.

To pass to the limit in this bound, we will use continuity of $G_j$. Since $F_j$ is continuous, $G_j$ is continuous on each open interval between consecutive integers. At a positive integer $k$, the left-hand limit is $k-1+F_j(1)=k$, while the value and the right-hand limit are both $k+F_j(0)=k$. Since $F_j(0)=0$ and $F_j$ is continuous at zero, $G_j$ is also continuous at zero. Thus $G_j$ is continuous on $[0,\infty)$.

When $j=1$, the conclusion of the following lemma follows from the exact result of Das, Gan and Sudakov~\cite{DasGanSudakov2016} for families of size $m=O(n)$. For a unified treatment of the critical scales, we nevertheless include a direct proof of this case.

\begin{lemma}\label{lem:all-critical}
Fix $1\le j\le r$ and $c>0$. For every sequence of integers $m=m(n)$ satisfying $0\le m\le\binom nr$ and $m/n^j\to c$ as $n\to\infty$, we have
\begin{equation}\label{eq:all-critical-limit}
M_r(n,m)\le\left(\frac{G_j(j!c)}{2(j-1)!}+o(1)\right)n^{j+1}.
\end{equation}
When $j=r$, we have $c\le1/r!$ because $m\le\binom nr$, and hence $G_r(r!c)=F_r(r!c)$.
\end{lemma}

\begin{proof}
When $j=r$, the assertion follows from Lemma~\ref{lem:dense} and the continuity of $F_r$, because $r!m/n^r\to r!c\le1$ and $G_r=F_r$ on $[0,1]$. We next consider $2\le j<r$. Fix $C>c$ and $0<\eta,\delta<1$, independently of $n$. Since $m/n^j\to c<C$, we have $m\le Cn^j$ for all sufficiently large $n$. For such $n$, take an arbitrary family $\mathcal F\subseteq\binom{[n]}r$ with $|\mathcal F|=m$. By Lemma~\ref{lem:all-extraction}, we can partition this family as $\mathcal F=\mathcal F_1\mathbin{\dot\cup}\cdots\mathbin{\dot\cup}\mathcal F_K\mathbin{\dot\cup}\mathcal F_0$. For every $1\le i\le K$, all members of $\mathcal F_i$ contain a fixed $(r-j)$-set $T_i$, and $K\le r\binom{r-1}{j-1}C/(\eta\delta)$. In particular, this upper bound on $K$ is independent of both $n$ and $\mathcal F$. Writing $m_i=|\mathcal F_i|$ for $0\le i\le K$, we have $\sum_{i=0}^K m_i=m$, and by~\eqref{eq:all-extraction},
\begin{equation}\label{eq:all-critical-extraction}
A(\mathcal F)\le\sum_{i=1}^K M_j(n,m_i)+\frac{j-1}{2}nm_0+O_{r,j}\!\left(mn\left(\eta+\frac{\delta^{1/(j-1)}}{\eta}\right)+Km\right).
\end{equation}

To estimate the terms $M_j(n,m_i)$, set $y_i=j!m_i/n^j$ for $0\le i\le K$. For $i\ge1$, deleting $T_i$ from each member of $\mathcal F_i$ produces a $j$-uniform family of size $m_i$ on $[n]\setminus T_i$, which can also be viewed as a family on $[n]$. Hence $m_i\le\binom nj$ and $0\le y_i\le1$; for the remaining part, we only require $y_0\ge0$. By Lemma~\ref{lem:dense}, for each $1\le i\le K$ we have
\[
M_j(n,m_i)=\frac{n^{j+1}}{2(j-1)!}F_j(y_i)+O_j(n^j),
\]
where the error is uniform in $m_i$. Recall that $\beta_j=(j-1)/j$. Since $y_0=j!m_0/n^j$, we have $(j-1)nm_0/2=n^{j+1}\beta_jy_0/(2(j-1)!)$. Substituting these expressions into~\eqref{eq:all-critical-extraction} and using $m\le Cn^j$, we obtain
\[
A(\mathcal F)\le\frac{n^{j+1}}{2(j-1)!}\left(\sum_{i=1}^K F_j(y_i)+\beta_jy_0\right)+O_{r,j,C}\!\left(n^{j+1}\left(\eta+\frac{\delta^{1/(j-1)}}{\eta}\right)\right)+O_{r,j,C,\eta,\delta}(n^j).
\]
The last error term includes both the sum $O_j(Kn^j)$ of the $K$ errors from Lemma~\ref{lem:dense} and the error $O_{r,j}(Km)$ for pairs in different parts. Both have the stated order because $K\le r\binom{r-1}{j-1}C/(\eta\delta)$ and $m\le Cn^j$. These bounds are uniform over all families $\mathcal F$ of size $m$.

By Corollary~\ref{cor:all-packing}, we obtain that $\sum_{i=1}^K F_j(y_i)+\beta_jy_0\le G_j(j!m/n^j)$, since $\sum_{i=0}^K y_i=j!m/n^j$. We now take the maximum over all $\mathcal F\subseteq\binom{[n]}r$ with $|\mathcal F|=m$, divide by $n^{j+1}$, and let $n\to\infty$ while $C,\eta,\delta$ remain fixed. Since $G_j$ is continuous and $m/n^j\to c$, we obtain
\[
\limsup_{n\to\infty}\left(\frac{M_r(n,m)}{n^{j+1}}-\frac{G_j(j!c)}{2(j-1)!}\right)
\le C_{r,j,C}\left(\eta+\frac{\delta^{1/(j-1)}}{\eta}\right).
\]
The constant $C_{r,j,C}$ is independent of $\eta$ and $\delta$. Keeping $C$ fixed, we first let $\delta\downarrow0$ with $\eta$ fixed, and then let $\eta\downarrow0$. We obtain the required upper bound.

For $j=1<r$, fix $0<\eta<1$ independently of $n$ and take an arbitrary family $\mathcal F\subseteq\binom{[n]}r$ of size $m$. As in the proof of Lemma~\ref{lem:all-extraction}, let $T_1,\ldots,T_K$ be all $(r-1)$-sets with $d_{\mathcal F}(T_i)\ge\eta n$. Since $m/n\to c$, we have $K=O_{r,c,\eta}(1)$ by the same codegree count.

Fix this ordering of $T_1,\ldots,T_K$ and, for $1\le i\le K$, define
\[
\mathcal F_i=\{E\in\mathcal F:T_i\subseteq E\text{ and }T_\ell\nsubseteq E\text{ for every }1\le\ell<i\}.
\]
Let $\mathcal F_0=\{E\in\mathcal F:T_i\nsubseteq E\text{ for every }1\le i\le K\}$. Then we have $\mathcal F=\mathcal F_1\mathbin{\dot\cup}\cdots\mathbin{\dot\cup}\mathcal F_K\mathbin{\dot\cup}\mathcal F_0$. Put $m_i=|\mathcal F_i|$ for $0\le i\le K$. For $i\ge1$, every member of $\mathcal F_i$ has the form $T_i\cup\{x\}$ with $x\notin T_i$, so $m_i\le n-r+1\le n$. Two distinct members have different additional elements and therefore intersect in exactly $T_i$. Consequently, $A(\mathcal F_i)=\binom{m_i}2\le m_i^2/2$.

We next bound $A(\mathcal F_0)$. For every $T\in\binom{[n]}{r-1}$, either $T$ is one of $T_1,\ldots,T_K$, in which case $d_{\mathcal F_0}(T)=0$, or $d_{\mathcal F_0}(T)\le d_{\mathcal F}(T)<\eta n$. By the codegree identity~\eqref{eq:codegrees}, we therefore have
\[
2A(\mathcal F_0)=\sum_{T\in\binom{[n]}{r-1}}d_{\mathcal F_0}(T)\bigl(d_{\mathcal F_0}(T)-1\bigr)\le\eta n\sum_{T\in\binom{[n]}{r-1}}d_{\mathcal F_0}(T)=r\eta nm_0.
\]
Using the same counting argument for pairs in different parts as in Lemma~\ref{lem:all-extraction}, we bound their number by $rKm$. To estimate the total contribution from $\mathcal F_1,\ldots,\mathcal F_K$, set $y_i=m_i/n$ for $0\le i\le K$. Then $y_0\ge0$, $0\le y_i\le1$ for $i\ge1$, and $\sum_{i=0}^K y_i=m/n$. Recall that $F_1(u)=u^2$ and $\beta_1=0$. By Corollary~\ref{cor:all-packing}, we obtain
\[
\sum_{i=1}^K\left(\frac{m_i}{n}\right)^2=\sum_{i=1}^K F_1(y_i)+\beta_1y_0\le G_1\!\left(\sum_{i=0}^K y_i\right)=G_1(m/n).
\]
We now add the bounds for adjacent pairs within each $\mathcal F_i$, within $\mathcal F_0$, and between different subfamilies. Since $A(\mathcal F_i)\le m_i^2/2$ for $i\ge1$ and $A(\mathcal F_0)\le r\eta nm_0/2$, we obtain from the preceding inequality and $m_0\le m$ that
\[
A(\mathcal F)\le\frac12\sum_{i=1}^K m_i^2+\frac{r\eta nm_0}{2}+rKm\le\frac{n^2}{2}G_1(m/n)+\frac{r\eta nm}{2}+rKm.
\]
Dividing by $n^2$ and letting $n\to\infty$ and then $\eta\to0$, we obtain the desired upper bound at the linear critical scale.
\end{proof}

Lemma~\ref{lem:all-critical} bounds $M_r(n,m)$ when $m/n^j$ tends to a positive finite limit. We next prove an upper bound when $n^{j-1}=o(m)$ and $m=o(n^j)$. In this range, we count the members containing many $(r-1)$-sets of high codegree directly, without first partitioning the family.

\begin{lemma}\label{lem:all-between}
Fix $2\le j\le r$. If $n^{j-1}=o(m)$ and $m=o(n^j)$, then
\begin{equation}\label{eq:all-between}
M_r(n,m)\le\left(\frac{j-1}{2}+o(1)\right)mn.
\end{equation}
\end{lemma}

\begin{proof}
Fix a family $\mathcal F$ of size $m$ and let $\mathcal S=\{S\in\binom{[n]}{r-1}:d_{\mathcal F}(S)\ge\eta n\}$, where $\eta\in(0,1)$ will be specified below. Again $|\mathcal S|\le rm/(\eta n)$. Let $D$ be the number of members of $\mathcal F$ that contain at least $j$ members of $\mathcal S$. Consider the $(r-1)$-uniform hypergraph with vertex set $[r]$ and edge set $\{[r]\setminus\{i\}:1\le i\le j\}$. Every vertex lies in at least $j-1$ of its edges, so weights $1/(j-1)$ form a fractional cover. Every member counted by $D$ supplies an embedding of this fixed hypergraph into $\mathcal S$, and embeddings chosen for different members have different vertex images. Applying Lemma~\ref{lem:all-entropy}, we obtain
\begin{equation}\label{eq:all-between-exceptional}
D\le\bigl((r-1)!|\mathcal S|\bigr)^{j/(j-1)}\le C_r\left(\frac{m}{\eta n}\right)^{j/(j-1)}.
\end{equation}
As in the proof of Lemma~\ref{lem:all-extraction}, we have $\sum_{S\in\mathcal S}d_{\mathcal F}(S)\le(j-1)m+rD$ by counting pairs $(E,S)$ with $E\in\mathcal F$, $S\in\mathcal S$ and $S\subseteq E$. For every $(r-1)$-set $S\notin\mathcal S$, we have $d_{\mathcal F}(S)<\eta n$, while $d_{\mathcal F}(S)\le n$ for every $(r-1)$-set $S$. Hence
\[
2A(\mathcal F)\le n((j-1)m+rD)+r\eta nm.
\]
Put $u=m/n^j$ and $\eta=u^{1/(2j-1)}$, which lies in $(0,1)$ eventually. Substituting these values into~\eqref{eq:all-between-exceptional}, we have
\[
\frac Dm\le C_r u^{1/(j-1)}\eta^{-j/(j-1)}=C_r u^{1/(2j-1)}=C_r\eta.
\]
Since $u\to0$, we obtain the required upper bound uniformly over all families of size $m$.
\end{proof}

The remaining sparse range is exact. If $0\le m\le n-r+1$, take $m$ distinct extensions of one $(r-1)$-set. Every pair is adjacent, and the total number of unordered pairs is an upper bound for any family. Hence
\begin{equation}\label{eq:all-sparse}
M_r(n,m)=\binom m2\qquad(0\le m\le n-r+1).
\end{equation}

To compare~\eqref{eq:P} with the estimates at the different scales, we need the behavior of $F_j(p)$ near zero and a uniform bound on $G_j(x)-x$. For each fixed $j\ge2$, we have $2p-1+(1-p)^{1+1/j}=\beta_jp+O_j(p^2)$ as $p\downarrow0$ by Taylor's theorem. Since $p^{1+1/j}=o(p)$ and $\beta_j>0$, this expression exceeds $p^{1+1/j}$ for all sufficiently small positive $p$. Hence
\begin{equation}\label{eq:all-profile-expansions}
F_j(p)=\beta_jp+O_j(p^2)\qquad(p\downarrow0,\ j\ge2).
\end{equation}
For $j=1$, we retain the exact identity $F_1(p)=p^2$. Since $0\le F_j(u)\le u$ on $[0,1]$, we obtain from the definition of $G_j$, for every $j\ge1$ and $x\ge0$, that
\[
G_j(x)-x=F_j(\{x\})-\{x\}\in[-1,0].
\]
Thus $G_j(x)=x+O(1)$ uniformly for $x\ge0$. We also need a lower bound. Since $0\le\beta_j\le1$, we obtain from~\eqref{eq:all-profile-bounds} that
\[
G_j(x)\ge\lfloor x\rfloor+\beta_j\{x\}=\beta_jx+(1-\beta_j)\lfloor x\rfloor\ge\beta_jx.
\]

The single expression~\eqref{eq:P} assembles these regimes without choosing a scale. The function $P_r(n,m)$ is positive throughout the stated range. Indeed, since $\sum_{j=2}^r(j-1)mn/2=r(r-1)mn/4$ and $G_j(x)\ge\beta_jx$, we have $P_r(n,m)\ge n^2G_1(m/n)/2-m/2$. Writing $m=qn+s$ with integers $q\ge0$ and $0\le s<n$, we have the exact identity
\[
\frac{n^2}{2}G_1(m/n)-\frac m2=q\binom n2+\binom s2.
\]
For $n\ge r\ge2$ and $m\ge2$, either $q\ge1$, or $q=0$ and $s=m\ge2$. The right side is positive in either case.

Now we are ready to prove Theorem~\ref{thm:main}.

\begin{proof}[Proof of Theorem~\ref{thm:main}]
Fix $r\ge2$ and consider integers $2\le m=m(n)\le\binom nr$. We first prove $M_r(n,m)\le(1+o(1))P_r(n,m)$ in each of the following three cases. In each case, we expand $P_r(n,m)$ and compare it with an upper bound for $M_r(n,m)$. We will then construct families attaining the corresponding lower bounds and use a subsequence argument to deduce uniformity over $2\le m\le\binom nr$.

\noindent\textbf{Case 1.} $m=o(n)$. Since $m/n<1$ for all sufficiently large $n$ and $G_1(x)=F_1(x)=x^2$ for $0\le x<1$, the term with $j=1$ in the sum defining $P_r(n,m)$ is $n^2G_1(m/n)/2=m^2/2$. For each $2\le j\le r$, we also have $j!m/n^j\to0$, so $G_j(j!m/n^j)=F_j(j!m/n^j)$ for all sufficiently large $n$. Applying~\eqref{eq:all-profile-expansions} with $p=j!m/n^j$, we obtain
\[
\frac{n^{j+1}}{2(j-1)!}G_j\!\left(\frac{j!m}{n^j}\right)
=\frac{j-1}{2}mn+O_r\!\left(\frac{m^2}{n^{j-1}}\right).
\]
Since $r$ is fixed, the sum of the error terms is $O_r(m^2/n)$. Substituting these estimates into~\eqref{eq:P} and using $\sum_{j=2}^r(j-1)=r(r-1)/2$, we obtain
\[
\begin{aligned}
P_r(n,m)
&=\sum_{j=1}^r\frac{n^{j+1}}{2(j-1)!}G_j\!\left(\frac{j!m}{n^j}\right)-\frac{r(r-1)}4mn-\frac m2\\
&=\frac{m^2}{2}+\frac{r(r-1)}4mn-\frac{r(r-1)}4mn-\frac m2+O_r(m^2/n)\\
&=\binom m2+O_r(m^2/n).
\end{aligned}
\]
Since $\binom m2\ge m^2/4$ for $m\ge2$, we have $P_r(n,m)=(1+o(1))\binom m2$, including when $m$ stays bounded. Since $m=o(n)$ and $r$ is fixed, we have $m\le n-r+1$ for all sufficiently large $n$. By~\eqref{eq:all-sparse}, $M_r(n,m)=\binom m2$. Therefore,
\[
\frac{M_r(n,m)}{P_r(n,m)}=\frac{\binom m2}{(1+o(1))\binom m2}=1+o(1),
\]
which establishes the required upper bound in this case.

\noindent\textbf{Case 2.} $m/n^k\to c\in(0,\infty)$ for some $1\le k\le r$. For $j<k$, using $G_j(x)=x+O(1)$, we obtain
\[
\frac{n^{j+1}}{2(j-1)!}G_j\!\left(\frac{j!m}{n^j}\right)
=\frac j2mn+O_r(n^{j+1}).
\]
For $k<j\le r$, we have $j!m/n^j=j!(m/n^k)n^{k-j}\to0$, since $m/n^k\to c$ and $k-j<0$. Thus, for all sufficiently large $n$, we have $0\le j!m/n^j<1$ and hence $G_j(j!m/n^j)=F_j(j!m/n^j)$. Applying~\eqref{eq:all-profile-expansions} with $p=j!m/n^j$, we obtain, for $j>k$,
\[
\frac{n^{j+1}}{2(j-1)!}G_j\!\left(\frac{j!m}{n^j}\right)=\frac{j-1}{2}mn+O_r\!\left(\frac{m^2}{n^{j-1}}\right).
\]
The sum of these errors is $O_{r,c}(n^k)$, and $m/2=O_{r,c}(n^k)$. Since $\sum_{j<k}jmn/2+\sum_{j>k}(j-1)mn/2=r(r-1)mn/4$, these terms cancel $-r(r-1)mn/4$ in the definition~\eqref{eq:P} of $P_r(n,m)$. Consequently,
\begin{equation}\label{eq:all-P-critical}
P_r(n,m)=\frac{n^{k+1}}{2(k-1)!}G_k\!\left(\frac{k!m}{n^k}\right)+O_{r,c}(n^k).
\end{equation}
Since $G_k$ is continuous, we have $P_r(n,m)/n^{k+1}\to G_k(k!c)/(2(k-1)!)$. By Lemma~\ref{lem:all-critical}, we therefore obtain $\limsup_{n\to\infty}M_r(n,m)/P_r(n,m)\le1$ at every critical scale, including the dense scale $k=r$.

\noindent\textbf{Case 3.} $n^{k-1}=o(m)$ and $m=o(n^k)$ for some $2\le k\le r$. For $1\le j<k$, we have $j!m/n^j=j!(m/n^{k-1})n^{k-1-j}\to\infty$, since $m/n^{k-1}\to\infty$ and $k-1-j\ge0$. Using $G_j(x)=x+O(1)$, we obtain, for $j<k$,
\[
\frac{n^{j+1}}{2(j-1)!}G_j\!\left(\frac{j!m}{n^j}\right)=\frac j2mn+O_r(n^{j+1}).
\]
For $k\le j\le r$, we have $j!m/n^j=j!(m/n^k)n^{k-j}\to0$, since $m/n^k\to0$ and $k-j\le0$. Thus, for all sufficiently large $n$, we have $0\le j!m/n^j<1$ and hence $G_j(j!m/n^j)=F_j(j!m/n^j)$. Applying~\eqref{eq:all-profile-expansions} with $p=j!m/n^j$, we obtain, for $j\ge k$,
\[
\frac{n^{j+1}}{2(j-1)!}G_j\!\left(\frac{j!m}{n^j}\right)=\frac{j-1}{2}mn+O_r\!\left(\frac{m^2}{n^{j-1}}\right).
\]
We now substitute these two estimates into~\eqref{eq:P} and use the identity $\sum_{j=1}^{k-1}jmn/2+\sum_{j=k}^r(j-1)mn/2-r(r-1)mn/4=(k-1)mn/2$. For $j<k$, we have $n^{j+1}\le n^k$, whereas for $j\ge k$, we have $m^2/n^{j-1}\le m^2/n^{k-1}$. Since $r$ is fixed, the errors in the two estimates sum to $O_r(n^k+m^2/n^{k-1})$. Including the remaining term $-m/2$ from~\eqref{eq:P}, we obtain
\begin{equation}\label{eq:all-P-between}
P_r(n,m)=\frac{k-1}{2}mn+O_r\!\left(n^k+\frac{m^2}{n^{k-1}}+m\right).
\end{equation}
The error divided by $mn$ tends to zero, since $n^{k-1}/m\to0$, $m/n^k\to0$, and $1/n\to0$. By Lemma~\ref{lem:all-between}, we again obtain $\limsup_{n\to\infty}M_r(n,m)/P_r(n,m)\le1$.

We now prove the matching lower bounds. If $2\le m=o(n)$, then eventually $m\le n-r+1$, and the family in~\eqref{eq:all-sparse} has exactly $\binom m2$ adjacent pairs. Thus the preceding expansion of $P_r(n,m)$ implies $\liminf_{n\to\infty}M_r(n,m)/P_r(n,m)\ge1$. If $m/n^r\to c>0$, Lemma~\ref{lem:dense} provides a family of size $m$ with $\bigl(F_r(r!c)/(2(r-1)!)+o(1)\bigr)n^{r+1}$ adjacent pairs. Since $r!c\le1$ and $G_r=F_r$ on $[0,1]$, we obtain the same lower bound for the ratio from~\eqref{eq:all-P-critical} with $k=r$.

For all other scales, we use the following construction. Fix an integer $1\le d<r$, let $q=q(n)\ge0$ be an integer with $q=o(n)$, and let $p\in[0,1)$. Choose $B\subseteq[n]$ with $|B|=r-d-1$ and distinct elements $v_1,\ldots,v_{q+1}$ outside $B$. Set $W=[n]\setminus(B\cup\{v_1,\ldots,v_{q+1}\})$ and $w=|W|=n-r+d-q$. These choices are possible for all sufficiently large $n$, since $q=o(n)$; moreover, $w=(1+o(1))n$.

We first define a family $\mathcal H\subseteq\binom Wd$. If $p^{1+1/d}\ge2p-1+(1-p)^{1+1/d}$, choose $U\subseteq W$ with $|U|=\lfloor p^{1/d}w\rfloor$ and set $\mathcal H=\binom Ud$. Otherwise, choose $V\subseteq W$ with $|V|=\lfloor(1-p)^{1/d}w\rfloor$ and set $\mathcal H=\binom Wd\setminus\binom Vd$, the family of all $d$-subsets of $W$ meeting $W\setminus V$. For $1\le i\le q$, define
\[
\mathcal F_i=\left\{B\cup\{v_i\}\cup E:E\in\binom Wd\right\}.
\]
For the last subfamily, define
\[
\mathcal F_{q+1}=\{B\cup\{v_{q+1}\}\cup E:E\in\mathcal H\}.
\]
Every member of $\mathcal F_i$ contains $v_i$ and no other element of $\{v_1,\ldots,v_{q+1}\}$, so these subfamilies are pairwise disjoint. Let $\mathcal C=\mathcal F_1\mathbin{\dot\cup}\cdots\mathbin{\dot\cup}\mathcal F_{q+1}$. Since $|B|+1+d=r$, this is an $r$-uniform family on $[n]$.

We next count its members and adjacent pairs. For each $i$, the map $E\mapsto B\cup\{v_i\}\cup E$ is injective, and hence
\begin{equation}\label{eq:all-construction-size}
s:=|\mathcal C|=q\binom wd+|\mathcal H|.
\end{equation}
For distinct $E,E'\subseteq W$ of size $d$, we have $|(B\cup\{v_i\}\cup E)\cap(B\cup\{v_i\}\cup E')|=r-d+|E\cap E'|$. Thus the two resulting $r$-sets are adjacent exactly when $|E\cap E'|=d-1$. Each member of $\binom Wd$ has $d(w-d)$ adjacent members in that family, so $A(\mathcal F_i)=d(w-d)\binom wd/2$ for $1\le i\le q$, and $A(\mathcal F_{q+1})=A(\mathcal H)$. Counting only pairs within the same subfamily, we obtain
\begin{equation}\label{eq:all-construction-count}
A(\mathcal C)\ge\sum_{i=1}^{q+1}A(\mathcal F_i)=\frac{qd(w-d)}2\binom wd+A(\mathcal H).
\end{equation}

It remains to calculate $|\mathcal H|$ and $A(\mathcal H)$ in its two defining cases. In the first case, $|U|=p^{1/d}w+O(1)$, so
\[
|\mathcal H|=\binom{|U|}{d}=\frac{|U|^d}{d!}+O_d(w^{d-1})=\frac{p}{d!}w^d+O_d(w^{d-1}).
\]
For each member of $\mathcal H$, an adjacent member is obtained by deleting one of its $d$ elements and adding one of the $|U|-d$ other elements of $U$. Each unordered adjacent pair is counted twice, and therefore
\[
A(\mathcal H)=\frac{d(|U|-d)}2\binom{|U|}{d}=\frac{p^{1+1/d}}{2(d-1)!}w^{d+1}+O_d(w^d).
\]
This identity also holds when $|U|<d$, since then $\mathcal H$ is empty.

In the second case, $|V|=(1-p)^{1/d}w+O(1)$, and hence
\[
|\mathcal H|=\binom wd-\binom{|V|}{d}=\frac{w^d-|V|^d}{d!}+O_d(w^{d-1})=\frac{p}{d!}w^d+O_d(w^{d-1}).
\]
Every member of $\mathcal H$ has exactly $d(w-d)$ adjacent members in $\binom Wd$. Summing these numbers over $\mathcal H$ counts each adjacent pair within $\mathcal H$ twice and each adjacent pair between $\mathcal H$ and $\binom Vd$ once. To count the latter pairs, start from a member of $\binom Vd$, delete one of its $d$ elements, and add an element of $W\setminus V$. There are exactly $d(w-|V|)\binom{|V|}{d}$ such pairs, each counted once. Consequently,
\[
\begin{aligned}
A(\mathcal H)
&=\frac{d(w-d)}2\left(\binom wd-\binom{|V|}{d}\right)-\frac{d(w-|V|)}2\binom{|V|}{d}\\
&=\frac{2p-1+(1-p)^{1+1/d}}{2(d-1)!}w^{d+1}+O_d(w^d).
\end{aligned}
\]
By the definition~\eqref{eq:F} of $F_d(p)$, we have therefore constructed $\mathcal H$ satisfying
\begin{equation}\label{eq:all-partial-construction}
|\mathcal H|=\frac{p}{d!}w^d+O_d(w^{d-1}),
\qquad
A(\mathcal H)=\frac{w^{d+1}}{2(d-1)!}F_d(p)+O_d(w^d).
\end{equation}
When $d=1$, both expressions defining $F_d(p)$ equal $p^2$, so the first case defines $\mathcal H$ to consist of $\lfloor pw\rfloor$ singletons. When $p=0$, the first case defines $\mathcal H$ to be empty. Thus no additional construction is needed in either case.

We now choose the parameters for a critical scale $m/n^k\to c>0$ with $1\le k<r$. Set $d=k$, $q=\lfloor k!c\rfloor$ and $p=k!c-q$. Here $q$ is independent of $n$, so $w=n-O_{r,c}(1)$. By~\eqref{eq:all-construction-size} and~\eqref{eq:all-partial-construction}, we have
\[
s=\frac{q+p}{k!}w^k+O_{r,c}(w^{k-1})=cn^k+O_{r,c}(n^{k-1})=(c+o(1))n^k.
\]
Since $q+F_k(p)=G_k(k!c)$, we also obtain from~\eqref{eq:all-construction-count} and~\eqref{eq:all-partial-construction} that
\[
A(\mathcal C)\ge\frac{q+F_k(p)}{2(k-1)!}w^{k+1}-O_{r,c}(w^k)=\left(\frac{G_k(k!c)}{2(k-1)!}+o(1)\right)n^{k+1}.
\]
We have $|s-m|=o(n^k)$ because both $s/n^k$ and $m/n^k$ tend to $c$. If $s<m$, add $m-s$ missing $r$-sets; if $s>m$, delete $s-m$ members. Each changed member belongs to at most $r(n-r)$ adjacent pairs in $\binom{[n]}r$, so these changes decrease the adjacency count by at most $r(n-r)|s-m|=o(n^{k+1})$. The resulting family has exactly $m$ members and retains the stated lower bound. Comparing with~\eqref{eq:all-P-critical}, we conclude that $\liminf_{n\to\infty}M_r(n,m)/P_r(n,m)\ge1$ at every such critical scale.

Finally, suppose $n^{k-1}=o(m)$ and $m=o(n^k)$ for some $2\le k\le r$. In the same construction, set $d=k-1$, $q=\lfloor d!m/n^d\rfloor$ and $p=0$. Then $\mathcal H=\varnothing$, $q\to\infty$ and $q=o(n)$, so $w/n\to1$. By~\eqref{eq:all-construction-size}, we have $s=q\binom wd\le qn^d/d!\le m$. Moreover,
\[
\frac{s}{m}=\frac{q}{d!m/n^d}\prod_{a=0}^{d-1}\frac{w-a}{n}=1+o(1),
\]
because $d!m/n^d\to\infty$ and $(w-a)/n\to1$ for every $0\le a<d$. Since $\mathcal H$ is empty, we obtain from~\eqref{eq:all-construction-count} that
\[
A(\mathcal C)\ge\frac{d(w-d)}2s=\left(\frac{k-1}{2}+o(1)\right)mn.
\]
Add $m-s$ missing $r$-sets to obtain a family of size $m$. Adding members cannot decrease the number of adjacent pairs, so~\eqref{eq:all-P-between} implies $\liminf_{n\to\infty}M_r(n,m)/P_r(n,m)\ge1$ between consecutive critical scales as well. Together with the upper bounds, these constructions prove $M_r(n,m)/P_r(n,m)\to1$ in each of the three types of sequences considered above.

For completeness, we deduce the uniform assertion from these sequence-wise conclusions, without assuming that an arbitrary cardinality sequence has one predetermined growth scale. Suppose that the uniform relative error did not tend to zero. There would then be $\varepsilon>0$, a sequence of integers $n\to\infty$, and integers $2\le m=m(n)\le\binom nr$ with $|M_r(n,m)/P_r(n,m)-1|\ge\varepsilon$. By successively selecting a subsequence for each $j=1,\ldots,r$, we may assume that every ratio $m/n^j$ converges in the compactified interval $[0,\infty]$. Let these limits be $\lambda_j$. Since $m\le\binom nr$, we have $\lambda_r\le1/r!$. If $\lambda_1=0$, the sparse argument applies. If some $\lambda_k$ is finite and positive, the critical-scale argument applies. In the only remaining case, every $\lambda_j$ is either zero or infinity, $\lambda_1=\infty$, and $\lambda_r=0$. Let $k$ be the first index with $\lambda_k=0$. Then $2\le k\le r$, $m/n^{k-1}\to\infty$, and $m/n^k\to0$, so the intermediate-scale argument applies. Each case contradicts the fixed discrepancy $\varepsilon$. Thus
\[
\sup_{2\le m\le\binom nr}\left|\frac{M_r(n,m)}{P_r(n,m)}-1\right|\longrightarrow0.
\]
The values at $m=0,1$ are zero. The constructions above, together with those in Lemma~\ref{lem:dense} and~\eqref{eq:all-sparse}, have exactly $m$ members after the stated adjustments. This completes the proof of Theorem~\ref{thm:main}.
\end{proof}

\section*{Declaration of generative AI use}

The author used OpenAI's GPT-5.6 Sol to explore proof strategies and assist with language polishing and proofreading. In particular, while the author was working on the inequality for the functions $h_d$ in Lemma~\ref{lem:moment}, GPT-5.6 Sol proposed a generalization to a broader class of functions, now stated in Lemma~\ref{lem:general-moment}, and provided a draft proof. The author subsequently verified and refined that proof. During the final revision, OpenAI's GPT-6 Astra was used for proofreading and language polishing. All mathematical results, arguments and proofs were independently reviewed and verified by the authors, who take full responsibility for the accuracy and content of this work.

\bibliographystyle{plain}
\bibliography{kw_references}

\end{document}